\documentclass[11pt]{article}
\usepackage{amssymb}
\usepackage{amsmath}
\usepackage[dvips]{epsfig}
\usepackage[small]{caption}
\usepackage{soul}
\usepackage{graphicx}
\usepackage{mathrsfs}
\usepackage{wrapfig}
\usepackage{caption}
\usepackage{subcaption}
\usepackage{float}
\usepackage[subnum]{cases}
\usepackage[colorlinks=true,linkcolor=blue,citecolor=red]{hyperref}

\newtheorem{Lemma}{Lemma}[section]
\newtheorem{Theorem}{Theorem}
\newtheorem{Proposition}[Lemma]{Proposition}

\newtheorem{Remark}[Lemma]{Remark}
\newtheorem{Definition}[Lemma]{Definition}

\newenvironment{Proof}[1][Proof]%
 {\begin{trivlist} \item[]{\textbf{#1.} }}%
{\hspace*{\fill}$\rule{.4\baselineskip}{.4\baselineskip}$\end{trivlist}}

 {\begin{trivlist}\item[]\textbf{Acknowledgments }}{\end{trivlist}}

\makeatletter\@addtoreset{figure}{section}\makeatother

\makeatletter \@addtoreset{equation}{section} \makeatother

\newcommand{\R}{\mathbb{R}}
\newcommand{\C}{\mathbb{C}}

\newcommand{\Z}{\mathbb{Z}}

\def\Re{\mathop{\mathrm{Re}}}

\newcommand{\bbF}{\mathbb{F}}
\newcommand{\bbL}{\mathbb{L}}

\newcommand{\tprime}{{\prime\prime\prime}}
\newcommand{\pprime}{{\prime\prime}}
\newcommand{\caL}{\mathcal{L}}
\newcommand{\caO}{\mathcal{O}}

\newcommand{\rmd}{\mathrm{d}}
\newcommand{\rme}{\mathrm{e}}
\newcommand{\rmi}{\mathrm{i}}
\newcommand{\id}{\mathrm{\,Id}\,}

\renewcommand{\leq}{\leqslant}
\renewcommand{\geq}{\geqslant}
\def\beq{\begin{equation}}
\def\eeq{\end{equation}}
\def\cF{{\cal F}}
\def\cG{{\mathcal G}}
\def\cR{{\mathcal R}}
\def\mrF{{\mathrm F}}
\def\veps{\varepsilon}
\def\oK{\overline{K}}

\newcommand{\Rmnum}[1]{\expandafter\@slowromancap\romannumeral #1@}

\newfam\bifam
\font\tenbi=cmmib10 scaled \magstep1 \font\sevenbi=cmmib10 at 11pt
\font\fivebi=cmmib10 at 6pt \textfont\bifam = \tenbi
\scriptfont\bifam = \sevenbi \scriptscriptfont\bifam= \fivebi

\begin{document}

\thispagestyle{empty}

\title{\bf Undulated Bilayer Interfaces in the Planar Functionalized Cahn-Hilliard Equation}

\author{Keith Promislow $^\dag $ and Qiliang Wu $^\ddag$\\
\textit{\small $^\dag$ Department of Mathematics, Michigan State University, 619 Red Cedar Road East Lansing, MI 48824}\\
\textit{\small $^\ddag$ Department of Mathematics, Ohio University, Morton 321, 1 Ohio University, Athens, OH 45701 }}

\date{\small \today}

\maketitle

\begin{abstract}
\noindent
Experiments with diblock co-polymer melts display undulated bilayers that emanate from defects such as triple junctions and endcaps, \cite{batesjain_2004}. Undulated bilayers are characterized by oscillatory perturbations of the bilayer width, which decay on a spatial length scale that is long compared to the bilayer width. We mimic defects within the functionalized Cahn-Hillard free energy by introducing spatially localized inhomogeneities within its parameters. For length parameter $\veps\ll1$, we show that this induces undulated bilayer solutions whose width perturbations decay on an $O\!\left(\veps^{-1/2}\right)$ inner length scale  that is long in comparison to the $O(1)$ scale that characterizes the bilayer width.
\end{abstract}


 \hrule
 {\small
%
%
%
  }


\section{Defect structures in amphiphilic morphology}
The Functionalized Cahn Hilliard free energy models the interaction of amphiphilic molecules with solvent. It traces its origins back to bilinear models of microemulsions of oil, water, surfactant derived from scattering data by Tubner and Strey, \cite{TS-87}. Their model represents the free energy through the density $u$ of the surfactant phase, and incorporates forth derivatives of this variable to match the decay of the scattering intensity with respect to wave number. Crucially, the coefficient of the second derivative term is negative, while the zeroth-order derivative term is positive.  Quadratic energies have linear variations, hence this model captures the linear response of the system to departures from spatial uniformity.  Surfactant systems are strongly phase separated, and not spatially homogeneous.  Gompper and Schick, \cite{GS-90} and later Gompper and Goos, \cite{GG-94} proposed nonlinear extensions of the energy that modulate the coefficients depending upon the local density, with the sign of the quadratic coefficient positive at very large or small densities of surfactant and negative a intermediate densities. This generically makes this coefficient a good candidate for the second derivative of a double-well potential.  Indeed, considering a domain $\Omega\subset{\mathbb R}^3$ subject to zero-flux or periodic boundary conditions, the energy introduced by Gompper and Goos can be written as a perturbation of a completed square in the highest order derivatives,
\beq
\cF_{GG}(u) :=\int_\Omega \frac{1}{2}\left(\veps^2 \Delta u-W'(u)\right)^2 +f_1 W(u)\,\mathrm d x.
\eeq
Here $0<\veps\ll1$ is a ratio of molecular and domain length scales, $W$ is a smooth double well potential with two local minima and $f_1$ is a bifurcation parameter that can take small positive or negative values. The quadratic term is often referred to as the Willmore energy, since for codimension one interfaces it generically reduces to the surface integral of the square of mean curvature. 

In \cite{PromislowWetton_2009} the model of Gompper-Goos was scaled and generalized into the functionalized Cahn-Hillard (FCH) energy
\beq
\label{e:FCHE}
\cF_{\mathrm FCH}(u):= \int_\Omega \frac12 \left( \veps^2 \Delta u - W'(u)\right)^2 - \veps^p\left(\eta_1 \frac{\veps^2}{2} |\nabla u|^2 + \eta_2 W(u)\right)\,
\mathrm d x.
\eeq
Here $\veps\ll1$ characterizes the ratio of the characteristic length of the surfactant molecules to the domain size and $W$ is an double well with unequal minima at $u=0$ and $u=u_+>u^\circ$ and a local maximum at $u=u^\circ.$ These satisfy $W(0)=0>W(u_+)$, and are non-degenerate in the sense that $W^{\prime\prime}(0)>0$,  $W^{\prime\prime}(u_+)>0$, and $W^{\prime\prime}(u^\circ)<0$. 
This form emphasizes the nearly ``perfect-square'' structure of the energy corresponding to $\eta_1=\eta_2=0,$ or equivalently to $p\to\infty.$ Indeed the value of $p$ represents a distinguished limit that slaves the Gompper's bifurcation parameter $f_1=-\veps^p\eta_2$ to $\veps$. The parameter $\eta_1>0$ incorporates the strength of the hydrophilicity of the solvent head groups. The choice $p=2$ yields an asymptotic balance between the functionalization terms, controlled by $\eta_1$ and $\eta_2$, and the residual of the dominant Willmore term. For $p=1$ these terms dominate the residual of the Willmore term. We remark that the $\veps^p\eta_1\veps^2|\nabla u|^2$ functionalization term can  be replaced with a $\veps^p\eta_1 u W'(u)$ potential, up to terms $O(\veps^{2p})$ by redefining the well shape $W$ at $O(\veps^p)$.

\begin{figure}[H]
 \vspace{-0.1in}
  \begin{center}
    \includegraphics[height=1.8in,angle=-0]{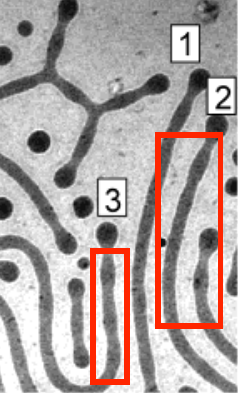}
  \end{center}
\vspace{-0.15in}
  \caption{Cryo-TEM images of blends of amphiphilic diblock polymer in water. 
  A mixture of diblocks with hydrophobic/hydrophillic chain lengths of 170/110 and 46/58, respectively. This mixture 
  produces visible undulations behind the endcaps, see the red boxes outlining the structures labeled '2' and '3'. Reprinted with permission from Figures 7\&8 of \cite{batesjain_2004}.}
\label{f:capy}
 \vspace{-5pt}
\end{figure}

The FCH free energy supports spatially extended structures that have a dimensional reduction. Posed in $\mathbb R^3$ these include codimension one bilayers and codimension two filaments; the codimension three micelles are not spatially extended in $\mathbb R^3$.  A defining feature of the nanoscale morphology produced by amphiphilic diblock polymers is a tendency for these spatially extended structures to undulate in the neighborhood of defects. Defects are defined as localized structures that break the dimensional reduction and include endcaps that terminate filaments or bilayers and $Y$ junctures. Undulations are long-range modulations of the thickness of the extended structure. Experiments reported in \cite{batesjain_2003, batesjain_2004} show that the modulations have wave-lengths that are comparable to the thickness of the structure and have amplitudes that attenuate on a length scale that is long compared to the interfacial thickness.  Undulations can be seen in 
in Figure\,\ref{f:capy} behind the end-cap defects that terminate the filament morphologies, particularly in the region within the red boxes near the endcaps labeled `2' and `3'. 

In this work we argue that breaking the perfect-square formulation of FCH free energy provides a mechanism that triggers the onset of undulations that decay on a slow, $O(\sqrt{\veps})$, length scale in the tangential direction.  We focus on flat bilayer critical points of the FCH in $\mathbb R^2$. In the spatial dynamics formulation we identify two pairs of purely imaginary eigenvalues within the linearization of the FCH in the perfect-square form. We show that these eigenvalues merge and either remain purely imaginary or bifurcate into four complex-conjugate eigenvalues with $O(\sqrt{\veps})$ real part when the perfect-square structure is broken at $O(\veps)$, see Figure\,\ref{f:spec}.
The former case, when the eigenvalues remain purely imaginary, was addressed in \cite{PW_14}, and corresponds to the creation of a family of bilayer profiles whose width is modulated by spatially periodic pearled patterns.  In the current work we address the the latter case, showing that the complex eigenvalues lead to the formation of undulations that form in presence of localized defects. Specifically we induce the defects by inserting spatial inhomogenities in $\eta_2=\eta_2(x)$. We establish that the bilayer solution persists under this perturbation, leading to an undulated equilibrium characterized by quasi-periodic oscillations with an $O(\veps)$ wavelength whose amplitude decays on an longer $O(\sqrt{\veps})$ length scale.


\paragraph{Bilayer solutions of the functionalized Cahn-Hilliard free energy.} 
We study the strong regime of the FCH, \eqref{e:FCHE}, with $p=1$ and subject to zero-flux boundary conditions and a mass constraint
$$\int_\Omega u\, dx = M.$$ 
 The critical points satisfy the Euler-Lagrange equation
 \begin{equation}
 \label{e:sFCH}
 (\veps^2\Delta - W^{\prime\prime}(u) + \veps \eta_1)
 (\veps^2\Delta u - W^{\prime}(u))+\veps\eta_d W^\prime(u)=\veps \gamma,
 \end{equation}
 where $\eta_d:=\eta_1-\eta_2$ and $\gamma\in\R$ is the Lagrange multiplier associated to conservation of mass of the FCH equation. 
 We assume the well $W$ is smooth and simplify the system, moving it to the plane, 
 \[
 \Omega=\R^2,
 \]
 and fixing a flat interface $\Gamma_f:=\{(x_1, 0)\mid x\in\R\}$ so that we may rewrite the Euler-Lagrange equation \eqref{e:sFCH} in the  in-plane/scaled-normal coordinates $(\tau,r)=(x_1, x_2/\veps)$, for which it takes the form
\begin{equation}
 \label{e:2sFCH}
 \left(\partial_r^2-W^{\prime\prime}(u)+\veps^2\partial_\tau^2+\veps\eta_1\right)
 \left(\partial_r^2u-W^\prime(u)+\veps^2\partial_\tau^2u\right)+\veps\eta_d(\tau) W^\prime(u)=\veps\gamma.
\end{equation}
We view the PDE as an infinite-dimensional dynamical system with $\tau=x_1$ playing the role of the evolution variable. The defect is induced by a spatial variation which we take in the form $\eta_2=\eta_{20}+\delta \eta_{21}(\tau)$ through the in-plane variable. This is made explicit in \eqref{eq-eta2-def}.
 
 In this scaling, the $\veps=0$ version of \eqref{e:2sFCH} possesses bilayer solutions. Both the in-plane variable $\tau$ and the functionalization parameters with their spatial variation are eliminated from the problem. The bilayers are the solutions of the second-order ODE
\beq\label{e:lODE}
\partial_r^2u-W^\prime(u)=0,
\eeq
that are homoclinic to origin, $u=0$, corresponding to the left local minimum of the well $W.$ The existence follows from classical planar dynamical systems techniques. We denote the unique
(up to translations)  orbit homoclinic to the origin, by $u_h=u_h(r;0)$, with the zero denoting the $\veps=0$ reduction.
The linearization of the system \eqref{e:lODE} at $u_0(r):=u_h(r;0)$ yields a Sturm-Liouville operator on the real line, 
\beq\label{e:L0}
\caL_0:=\partial_r^2-W^{\pprime}(u_0): H^2(\R)\longrightarrow L^2(\R),
\eeq
whose spectrum, according to standard Sturm-Liouville theory, consists of a simple positive eigenvalue $\lambda_0$ and the simple eigenvalue $\lambda_1= 0$ with the remainder lying strictly on the negative real axis. We denote the associated $L^2(\R)$ normalized eigenfunctions as $\psi_0$ and $\psi_1$, respectively. We also note that, unless otherwise noted, the notation $L^2(\R)$ is reserved for function of $r$ and $\langle\cdot,\cdot \rangle_{L^2(\mathbb{R})}$ denotes an integration over $r$.

In the case of spatially homogeneous $\eta_1$ and $\eta_2$ one may drop the $\tau$ derivatives and study the persistence of these homoclinic orbits for $0<\veps\ll1$. The work \cite{doelmanpromislow_2014} considered a flat interface and constructs bilayer profiles which are homoclinic to the far-field value
\beq
\label{eq-uinf}
u_\infty = \veps \frac{\gamma}{(W''(0))^2} + O(\veps^2),
\eeq
which is the unique small solution of the far-field equation
\[
(W^{\prime\prime}(u_\infty)-\veps \eta_2)W^\prime(u_\infty)=\veps \gamma.
\]
Introducing the constant $\eta_{d,0}=\eta_{2,0}-\eta_1,$ they establish the persistence of bilayer profiles for $\veps$ sufficiently small.
\begin{Theorem}[\cite{doelmanpromislow_2014}-Theorem 3.1]
\label{d:bilayer}
Fix $\gamma_0>0$,  then there exists $\veps_0>0$ such that for all $|\gamma|<\gamma_0$ and all $\veps\in[0,\veps_0)$, the constant coefficient bilayer ODE
\beq\label{e:rODE}
 \left(\partial_r^2-W^{\prime\prime}(u)+\veps\eta_1\right)
 \left(\partial_r^2u-W^\prime(u)\right)+\veps\eta_{d,0} W^\prime(u)=\veps\gamma,
\eeq
admits, up to translation, a unique solution $u_h(r;\veps)$, called the bilayer solution to the FCH, that is homoclinic to $u_\infty(\veps;\gamma).$
\end{Theorem}


\paragraph{Undulated bilayer interfaces induced by amphiphilic inhomogeneity.}

The FCH parameter $\eta_2$ is well-known to tune the energetic preference of the system for various codimensional morphologies, \cite{CKP_2019} and is the central bifurcation parameter in the formulation in \cite{GG-94}.  Our central result is that if the key parameter $\alpha_0$ defined in \eqref{e:alpha0} is negative, then spatially inhomogeneity in the parameter $\eta_2$ will induce long oscillations characteristic of the structures observed experimentally behind endcap defects in Figure\,\ref{f:capy}. Specifically for $0<\delta\ll1$ we consider inhomogeneity's
\beq\label{eq-eta2-def}
\eta_2(\tau; \delta)=\eta_{2,0} + \delta\xi\left(\sqrt{\lambda_0}\frac{\tau}{\veps}\right),
\eeq
 where $\xi$ has compact variation: $\xi'$ is a smooth function satisfying
\beq
\label{def-xi}
\xi'(t)=0, \quad |t|>T.
\eeq
If $\xi'$ has zero mass, then $\xi$ has identical limits as $t\to\pm\infty$ and we say $\xi$ is a localized inhomogeneity. Conversely, if $\xi'$ has non-zero mass, then  $\xi(\pm \infty)$ differ and we say $\xi$ is a transitional inhomogeneity. 
The impact of $\xi$ on the perturbed solution to \eqref{e:2sFCH} is characterized by the Fourier coefficients
 \beq\label{e:xi-Fourier}
 \begin{aligned}
 \Xi_{\mathrm{o},1}&:= \int_{\R} \xi'(t)\cos(t)\rmd t, \\ 
 \Xi_{\rme,1}&:=-\int_\R \xi'(t)\sin(t)\rmd t. 
 \end{aligned}
 \eeq

We establish the continuous bifurcation of undulated bilayer interfaces out of bilayer solutions  for $\veps>0$ and $\delta>0$ sufficiently small. The main significance is the appearance of the slow $O(\veps^{-1/2})$ decay of the undulations induced by the localized inhomogeneity in $\eta_2.$

\begin{Theorem}
\label{thm-main}
Assume that the pearling bifurcation parameter $\alpha_0$, defined in \eqref{e:alpha0}, is negative; that is, $\alpha_0<0$.
If in addition the following generic conditions hold 
\begin{enumerate}
    \item The Fourier coefficients $\Xi_1:=(\Xi_{\rme, 1}, \Xi_{\mathrm{o},1})$ are not identically zero, $|\Xi_1|>0$;
    \item The scaled inner-product $\beta_0:=\langle \psi_0, W^\prime(u_0)\rangle_{L^2(\R)}$ is non-zero. Here $u_0$ is the bilayer solution of \eqref{e:lODE} and $(\lambda_0,\psi_0)$ is the ground state eigen-pair of the associated linearization ${\cal L}_0$,  defined in \eqref{e:L0}. 
\end{enumerate}
Then for any $q>3/4$ there exist $\delta_0,\veps_0 >0$,  such that for all $(\delta\veps^{-q},\veps)\in(-\delta_0,\delta_0)\times(0,\veps_0)$, the stationary FCH equation \eqref{e:sFCH} with amphiphilic inhomogeneity
 $\eta_d(x_1,\delta)=\eta_{2,0}-\eta_1+\delta\xi(\sqrt{\lambda_0}\frac{x_1}{\veps})$, admits an undulated bilayer solution
\begin{equation}
\label{e:main-result}
\begin{aligned}
    u_n(x)= & u_{h}\left(\frac{x_2}{\veps};\veps,\eta_d(x_1,\delta)\right)+\frac{\beta_0 |\Xi_1|}{4\lambda_0^2\sqrt{-\alpha_0}}\delta\sqrt{\veps}e^{-\sqrt{-\alpha_0\lambda_0}\frac{x_1}{\sqrt{\veps}}}\cos\left(k(\veps)\frac{x_1}{\veps}+\Theta_1\right)\psi_0\left(\frac{x_2}{\veps}\right)+ \\
    & \hspace{0.5in} \caO(\delta\veps^{3/4}+\delta^2\veps^{-1/2}),
\end{aligned}
\end{equation}
 where $u_{h}$ is the $x_2$ dependent bilayer solution of Theorem \ref{d:bilayer} modulated by the $x_1$ variation in $\eta_d,$ see also equation \eqref{e:MBI}. The scaled wavenumber $k(\veps):=\sqrt{\lambda_0}A(\veps)=\sqrt{\lambda_0}+\caO(\veps)$, with $A$ defined in \eqref{e:AB}. The phase shift $\Theta_1$ is the angle of the vector $\Xi_1$. In addition, the error terms are taken in $H^4(\R^2)$ as functions of the inner coordinates $(t,r):=(\sqrt{\lambda_0}\frac{x_1}{\veps},\frac{x_2}{\veps})$.
\end{Theorem}

This result requires that $\beta_0\neq0.$
The homoclinic orbit $u_0$ solves \eqref{e:lODE} while $\psi_0>0$ is the ground-state eigenfunction of ${\cal L}_0$, defined in \eqref{e:L0}, corresponding to eigenvalue $\lambda_0>0.$ Moreover, the first excited-state eigenfunction $\psi_1=u_0'$ has eigenvalue $\lambda_1=0.$
Consequently we may write
\beq
\label{beta0-2}
\beta_0:= \langle\psi_0,W^\prime(u_0) \rangle_{L^2(\R)}=\langle\psi_0,\partial_r^2u_0 \rangle_{L^2(\R)}=-\langle\partial_r\psi_0,\psi_1 \rangle_{L^2(\R)}.
\eeq
The operator ${\cal L}_0$ is Sturmian, so by classical Sturm-Liouville theory all eigenvalues are simple and $\psi_0$ has even parity about $r=0$. Moreover $\psi_1$ has odd parity about $r=0$, with a simple zero at $r=0$ and is positive on $r>0$ and $r<0.$ If $\psi_0$ is monotonic on $r>0$ then we may deduce that $\beta_0\neq0.$  The results of \cite{PY-14} shown that as $W$ approaches an equal-depth well, then $u_0$ approaches a heteroclinic connection and $\psi_0\rightarrow |\psi_1|.$ In this limit we have $\beta_0\rightarrow 0,$ thus a non-zero value of $\beta_0$ is not immediate.  

The following result shows that for a significant class of wells $W$, the inner product $\beta_0$ is negative.
\begin{Lemma}
Assume that the ground state eigenvalue $\psi_0$ of ${\cal L}_0$ is scaled so that $\psi_0>0$. Let $u_{\rm max}\in(u^\circ,u_+)$ denote the smallest positive zero of $W$. If $W'''(u)<0$ for $u\in(0,u_{\rm max})$ then the inner product $\beta_0$ is strictly negative, in particular it is non-zero.
\end{Lemma}
\begin{Remark}
It is straight-forward to construct tilted double-well potentials $W$ satisfying the conditions imposed after \eqref{e:FCHE} for which $\beta_0$ defined in Lemma 1.1 is negative. Indeed the function
\[
W(u)=u^2(u-u_{\rm max})(u-cu_{\rm max}),
\]
does so for all $u_{\rm max}>0$ and all $c>3.$
\end{Remark}
\begin{Proof}
We take $\partial_r$ of the eigenvalue equation for $\psi_0$ and use $u_0'=\psi_1$ to obtain the identity,
\[ {\cal L}_0 \psi_0' = W'''(u_0)\psi_0\psi_1 +\lambda_0\psi_0'.\]
By the Fredholm alternative the right-hand side of this identity is orthogonal to $\psi_1$, which spans the kernel of ${\cal L}_0.$ Taking the inner product of the right-hand side with $\psi_1$, and using \eqref{beta0-2}  we find
\[ \beta_0 = \frac{1}{\lambda_0}\int_\R W'''(u_0)\psi_0\psi_1^2\, {\textrm d}r. 
\]
By assumption $W'''(u_0)<0$ while $\psi_0>0,$ and we conclude that $\beta_0\neq0.$ Since $\lambda_0>0$ we establish the result.
\end{Proof}



\begin{Remark}
The inhomogeneity $\xi'$ is chosen to have compact support as this allows for a simplification of the leading order terms in \eqref{e:main-result}. For $\xi'\in L^1(\R)$ a similar asymptotic form holds with an adjusted scaling with respect to $\veps$; see Lemma \ref{l:GK0} and the estimate \eqref{e:Gf-exact} in the proof of Theorem \ref{thm-main} for details.
\end{Remark}
\begin{Remark}
For $\alpha_0>0$ the unperturbed system supports pearled solutions that are perturbations of bilayers with spatially periodic variations in the bilayer width, see \cite{PW_14}. For the perturbed system the presence of the spatial inhomogeneity in $\eta_2$ generically excites resonant modes in the linear system that lead to secular growth of the underlying perturbation as measured in distance along the bilayer. Such growth is often saturated by the higher-order nonlinear terms. Consequently the inhomogenous system may support pearled patterns with defects,  but this analysis is outside the scope of our current framework.
\end{Remark}

\section{Center manifold reduction of bilayer profiles}
\label{s:CM-reduction}

For the flat interface $\Gamma$ and spatially constant parameters $\eta_1$ and $\eta_2$, the bilayer profiles constructed in Theorem\,\ref{d:bilayer} naturally extend to functions defined on the whole spatial domain. We call these functions bilayer interfaces, and their dynamic stability has been studied in \cite{doelmanpromislow_2014, PW_14}, which showed that they may be unstable to either pearling or meandering bifurcations depending upon parameter values. Pearling bifurcations correspond to high-frequency, periodic modulations of the through-plane on the fast $O(\veps)$ length scale. Meander bifurcations modulate the shape of the center line  bilayer interface, perturbing it from its the flat shape. These are generically long-wave effects with $O(1)$ spatial variation. 

A complete center manifold reduction that characterizes the possible pearled equilibrium local to the flat bilayer interface was developed in \cite{PW_14} via a spatial dynamics analysis.   We summarize these results as they provide a framework that motivates the genesis of the undulated bilayer interfaces constructed in section 3. Recalling $u_h$ constructed in Theorem\,\ref{d:bilayer},
we introduce the perturbation 
\[
v:=u-u_h,
\]
 and the linear operator 
 \beq
 \label{e:caLh-def}
 \mathcal{L}_h:=\partial_r^2-W^{''}(u_h).
 \eeq
 We change variables to $\tau=\veps t/\sqrt{\lambda_0}$, where $\lambda_0>0$ is the ground state eigenvalue of $\caL_0$. This is equivalent to 
 \beq \label{rescaled0}
(t,r)=(\sqrt{\lambda_0}x_1/\veps, x_2/\veps),
\eeq
for which \eqref{e:2sFCH} takes the form
\beq
\label{e:in-sFCH-2-0}
    (\partial_r^2-W^\pprime(u)+\lambda_0\partial_t^2+\veps\eta_1)(\partial_r^2u-W^\prime(u)+\lambda_0\partial_t^2u)+\veps\eta_{d,\delta}(t)W^\prime(u)=\veps\gamma.
\eeq
We now apply the spatial dynamics approach to rewrite \eqref{e:in-sFCH-2-0} as an infinite-dimensional dynamical system, where the rescaled in-plane variable $t$ is viewed as the evolution variable, bilayer solutions as equilibria and pearled bilayers as periodic temporal oscillations to these bilayer equilibria. More specifically, we denote $\cdot=\frac{\rmd}{\rmd t}$, introduce the variables
\begin{equation}
    V=\begin{pmatrix}
    v_1 \\ v_2 \\ v_3\\ v_4
    \end{pmatrix}:=\begin{pmatrix}v \\ \partial_tv \\  \mathcal{L}_hv+\lambda_0\partial_t^2v \\ \partial_t\big( \mathcal{L}_hv+\lambda_0\partial_t^2v \big)
    \end{pmatrix},
\end{equation}
and rewrite \eqref{e:in-sFCH-2-0} as an infinite-dimensional dynamical system 
 \begin{equation}\label{e:idds}
      \dot{V}=\bbL(\veps) V +\bbF(V,\veps),
 \end{equation}
 where the linear and strictly nonlinear terms are
 \beq\label{e:vLF}
  \bbL(\veps):=\begin{pmatrix} 0 & 1 & 0 & 0 \\ -\caL_h/\lambda_0 & 0 & 1/\lambda_0 & 0 \\ 0 & 0  & 0 & 1 \\ \mathcal{M}& 0 & -(\caL_h+\veps\eta_1)/\lambda_0 & 0  \end{pmatrix},\quad
  \bbF(V, \veps):=\begin{pmatrix}0 \\ 0\\ 0\\ \mathcal{F}(V, \veps)\end{pmatrix},
 \eeq
 the $(4,1)$ entry of $\bbL$ takes the form 
 \[
 \mathcal{M}:=-[\varepsilon\eta_dW^{\pprime}(u_h)-\left(\partial_r^2u_h-W^\prime(u_h)\right)W^{\tprime}(u_h)]/\lambda_0, \]
 and the nonlinearity $\mathcal{F}$ is given by
 \[
 \begin{aligned}
\mathcal{F}(V,\varepsilon):=&W^{\tprime}(u_h+v)\left(\partial_tv\right)^2+2\left(W^{\prime\prime}(u_h+v)-W^{\prime\prime}(u_h)\right)\partial_t^2v+\\
 &\left[\caL_h+\varepsilon\eta_2-\left(W^{\pprime}(u_h+v)-W^{\pprime}(u_h)\right)\right]\left(W^\prime(u_h+v)-W^\prime(u_h)-W^{\pprime}(u_h)v\right)/\lambda_0+\\
 &\left(W^{\pprime}(u_h+v)-W^{\pprime}(u_h)\right)\caL_h v-
 \left(\partial_r^2u_h-W^\prime(u_h)\right)\left(W^{\pprime}(u_h+v)-W^\pprime(u_h)-W^\tprime(u_h)v\right)/\lambda_0.
 \end{aligned}
 \]
 The spectrum of $\bbL(\veps)$ in \eqref{e:vLF} is determined from the eigenvalue problem
 \[
 \bbL (\veps)V=\lambda V.
 \]
 The operator $\bbL(\veps)$ is the vectorized version of the scalar operator
 \beq
 \label{Lelam-def}
 L(\veps,\lambda):= \left(\caL_h+\veps\eta_1+\lambda_0\lambda^2\right)\left(\caL_h+\lambda_0\lambda^2\right)+\veps(\eta_d W^{\pprime}(u_h)-W^\tprime(u_h)w_h),
 \eeq
 where we have introduced $w_h:=\left(\partial_r^2u_h-W^\prime(u_h)\right)/\veps$.  
 Accordingly the spectrum of $\bbL(\veps)$ agrees, up to multiplicity to the nontrivial solutions $\lambda$ of
 \beq\label{e:1eigen}
L(\veps,\lambda)v=0.
 \eeq
 Since $L(\cdot,\cdot)$ is connected to the Hessian of the FCH energy, it is natural that for $\veps=0$ it becomes a square of a second order operator, and the eigenvalue problem reduces to
 \[
 L(0,\lambda)v=(\caL_0+\lambda_0\lambda^2)^2v=0.
 \]
These observations imbue the spectrum of $\bbL(\veps)$ in $(L^2(\R))^4$, denoted $\sigma(\bbL(\veps))$, with the following properties:
\begin{enumerate}
\item $\sigma(\bbL(\veps))$ is symmetric with respect to the real and imaginary axis.
 \item $\sigma(\bbL(0))=\{ \pm\sqrt{-\lambda}\in\C\mid  \lambda\in\sigma(\caL_0/\lambda_0) \}=\{0, \pm\rmi\}\cup \{\pm\sqrt{-\lambda}\mid \lambda<0, \lambda\in\sigma(\caL_0/\lambda_0)\}$.
 \item The eigenvalue $\lambda=0$, called the meandering eigenvalue of $\bbL(0)$, has algebraic multiplicity $4$; $\lambda=\pm\rmi$ are the pearling eigenvalues of $\bbL(0)$, each has algebraic multiplicity $2$.
 \item  The eigenfunction $\psi_0$ and $\psi_1$ of $\caL_0$, defined in \eqref{e:L0}, satisfy $L(0,0)\psi_1=0$ and $L(0,\pm\rmi)\psi_0=0$.
\end{enumerate}
In the sequel we show that the continuation of the pearling eigenvalues $\pm\rmi$ as $\veps$ increases from zero determines much of the structure of the perturbed problem we study in section 3.  The double multiplicity of the pearling modes precludes a direct application of the implicit-function-theorem argument. A remedy, based on the observation that \eqref{e:1eigen} admits the expansion
\beq
\label{e:L-cal0}
L(\veps,\lambda)v=\left[ (\caL_0+\lambda_0\lambda^2)^2+\caO(\veps)\right]v,
\eeq
is to unfold the degeneracy through the change of variable
\[
\lambda^2=-1+\sqrt{\veps}\Lambda, \quad v=\psi_0+\veps \Psi, \quad \text{where}\,\langle \Psi, \psi_0\rangle_{L^2(\R)}=0.
\]
This allow us to recast \eqref{e:1eigen} as the search for the zeros of $F$ defined by
\beq\label{e:2eigen}
F(\Lambda, \Psi; \sqrt{\veps}):=\veps^{-1}L\left(\veps,\sqrt{ -1+\sqrt{\veps}\Lambda}\right)(\psi_0+\veps \Psi)=0.
\eeq
In the limit $\veps\to0^+$ it is straightforward to calculate that $F$ reduces to 
\beq\label{e:3eigen}
F(\Lambda, \Psi; 0)=(\caL_0-\lambda_0)^2\Psi+\lambda_0^2\Lambda^2\psi_0+(\caL_0-\lambda_0)(W^\tprime(u_0)u_1\psi_0)+(\eta_d W^\pprime(u_0)-W^\tprime(u_0)w_0)\psi_0,
\eeq
where we have introduced 
\beq\label{e:u1w0}
\begin{aligned}
u_1:=&\lim_{\veps\to0+}\frac{u_h-u_0}{\veps},\\
w_0:=&\lim_{\veps\to0+}w_h=\lim_{\veps\to0+}\frac{\partial_r^2u_h-W^\prime(u_h)}{\veps}=\caL_0^{-1}(\gamma-\eta_d W^\prime(u_0)).
\end{aligned}
\eeq

The quantity $\caL_0^{-1}(\gamma-\eta_d W^\prime(u_0))$ is well-defined since the operator $\caL_0$, is invertible on functions with even parity. Indeed expanding $u_h$
\begin{equation}\label{e:Taylor-uh}
    u_h=u_0+\veps u_1+\caO(\veps^2),
\end{equation}
in \eqref{e:rODE}, we deduce that
\begin{equation}
\label{e:u1}
\caL_0^2 u_1=\gamma-\eta_d W^\prime (u_0),
\end{equation}
and hence
\begin{equation}\label{e:w0u1}
    w_0=\caL_0u_1.
\end{equation}
The eigenvalue $\lambda_0$ of $\caL_0$ is geometrically simple with corresponding normalized eigenfunction $\psi_0$. The pearling parameter
\begin{equation}\label{e:alpha0}
    \alpha_0:=\frac{1}{4\lambda_0^2}\left\langle (W^\tprime(u_0)w_0-\eta_d W^\pprime(u_0))\psi_0, \psi_0 \right\rangle_{L^2(\R)}=\int_{\R}\left(W^\tprime(u_0)\caL_0u_1-\eta_dW^\pprime(u_0)\right)\psi_0^2\rmd r,
\end{equation}
is a zero of $F$, satisfying
\[
F(\pm2\sqrt{\alpha_0},\Psi_0; 0)=0,
\]
where 
$$\Psi_0:=(\caL_0-\lambda_0)^{-2}\Big[-4\lambda_0^2\alpha_0\psi_0-(\caL_0-\lambda_0)(W^\tprime(u_0)u_1\psi_0)+(W^\tprime(u_0)\caL_0u_1-\eta_d W^\pprime(u_0))\psi_0\Big],$$
see \cite{PW_14} for details.
\begin{Remark}
From \eqref{e:u1} the pearling parameter $\alpha_0$ can be written as 
\[\alpha_0=  \alpha_{01} \gamma  - \alpha_{02} \eta_d,
\]
where
\begin{equation}
\begin{aligned}
  \alpha_{01} &:= \frac{1}{4\lambda_0^2}\int_{\R} W^\tprime(u_0)(\mathcal{L}_0^{-1} 1)\psi_0^2\rmd r, \\
\alpha_{02} & :=\frac{1}{4\lambda_0^2}\int_{\R} \left(\mathcal{L}_0^{-1}W'(u_0)+  W^\pprime(u_0)\right)\psi_0^2\rmd r.
\end{aligned}
\end{equation}
The constants $\alpha_{01}$ and $\alpha_{02}$ depend only upon the form of the double well potential, $W$.
\end{Remark}
Moreover, the derivative $\nabla_{\Lambda,\Psi} F(\pm2\sqrt{\alpha_0},\Psi_0;0)$ is bounded and invertible, and for $0<\veps\ll1$ the implicit function theorem shows that \eqref{e:2eigen} admits solutions $(\Lambda_\pm,\Psi_\pm)$ with following expansions
\[
\Psi_\pm(r; \sqrt{\veps})=\Psi_0(r)+\caO(\veps), \quad \Lambda_\pm(\sqrt{\veps})=\pm2\sqrt{\alpha_0}+\caO(\sqrt{\veps}).
\]
These results are a reformulation of Lemma 2.9 in \cite{PW_14} and summarized in the following proposition.

\begin{figure}
    \centering
    \includegraphics[width=0.5\textwidth]{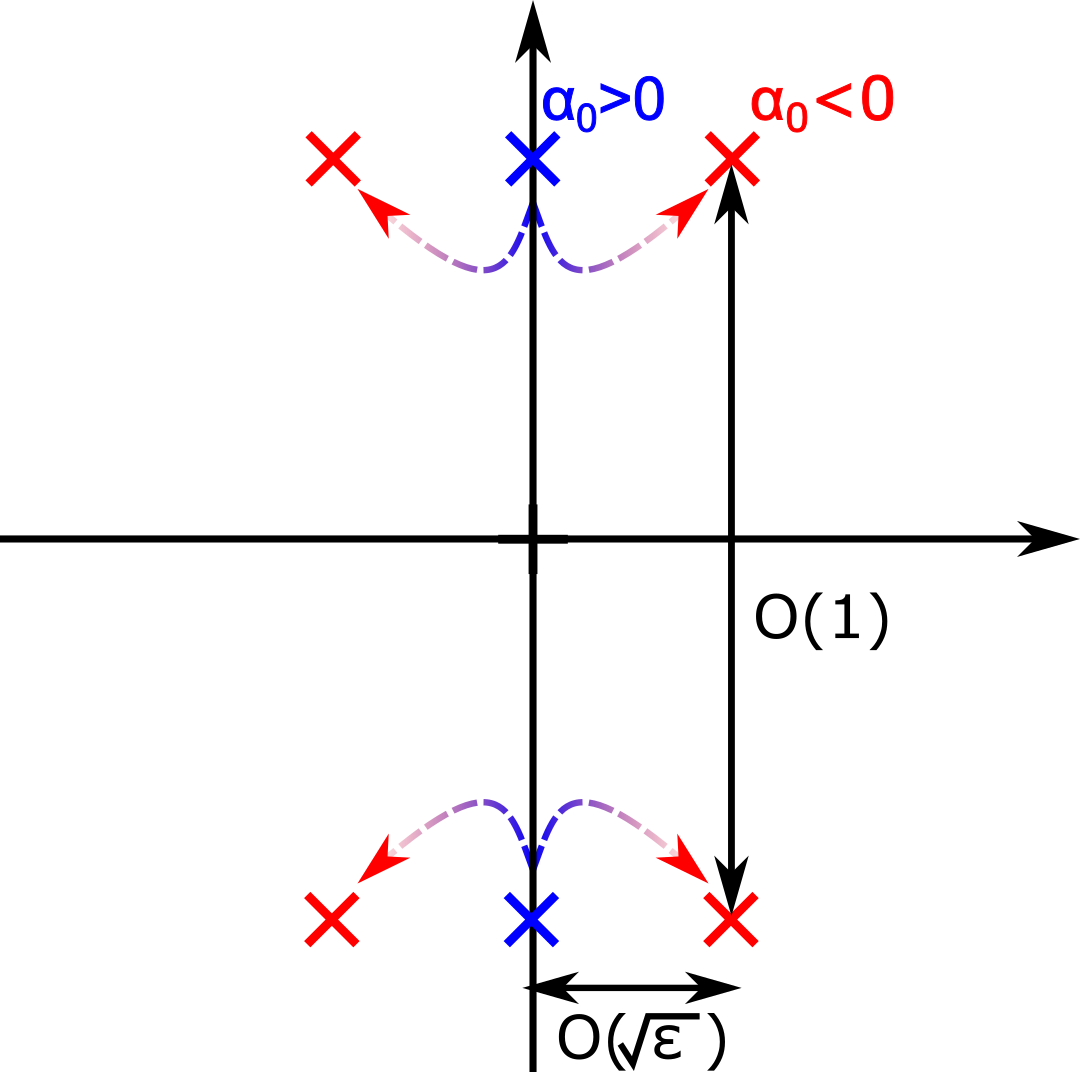}
    \caption{The operator $\bbL(\veps)$ admits two purely imaginary eigenvalues (blue crosses) with algebraic multiplicity 2 when $\veps=0.$ Given $\veps>0$, for $\alpha_0>0$ they remain purely imaginary while for $\alpha_0<0$ they split into four geometrically simple modes (red crosses).}
    \label{f:spec}
\end{figure}

\begin{Proposition} 
For $\veps>0$ sufficiently small, the operator $\bbL(\veps)$ admits four eigenvalues $\pm\lambda_p(\sqrt{\veps})$, $\pm\overline{\lambda_p}(\sqrt{\veps})$ with the expansion 
\[
\lambda_p(\sqrt{\veps})=\sqrt{-1+2\sqrt{\alpha_0}\sqrt{\veps}+\caO(\veps)}=\rmi+\sqrt{-\alpha_0}\sqrt{\veps}+\caO(\veps).
\]
The corresponding eigenfunction with respect to $\lambda_p$, denoted as $\Psi_p$, takes the form
\[
\Psi_p(r;\sqrt{\veps})=\psi_0(r)+\veps\Psi_0(r)+\caO(\veps^{3/2}),
\]
where we recall that $\psi_0$ is the normalized eigenfunction of $\caL_0$ with respect to $\lambda_0$ and 
\[
\Psi_0:=(\caL_0-\lambda_0)^{-2}\Big[-4\lambda_0^2\alpha_0\psi_0-(\caL_0-\lambda_0)(W^\tprime(u_0)u_1\psi_0)+(W^\tprime(u_0)\caL_0u_1-\eta_d W^\pprime(u_0))\psi_0\Big].
\]
Moreover, we have the following distinctive scenarios.
\begin{enumerate}
\item (Pearling) If $\alpha_0>0$, then the four eigenvalues  $\pm\lambda_p$, $\pm\overline{\lambda_p}$ are pure imaginary, giving rise to pearling bifurcation; see \cite{doelmanpromislow_2014, PW_14} and Fig. \ref{f:spec}.
\item (Undulations) If $\alpha_0<0$, then $\Re\lambda_p=\sqrt{-\alpha_0}\sqrt{\veps}+\caO(\veps)>0$ and the four eigenvalues $\pm\lambda_p$, $\pm\overline{\lambda_p}$ are geometrically simple. The eigenvalues $\lambda_p$ and $\overline{\lambda_p}$ are the leading modes of the unstable spectra of $\bbL(\veps)$ while the eigenvalues $-\lambda_p$ and $-\overline{\lambda_p}$ are the leading modes of the stable spectra of $\bbL(\veps)$; see Fig. \ref{f:spec} for an illustration.
\end{enumerate}
\end{Proposition}

In the case $\alpha_0<0$ we will show that, in the presence of defects, the eigenspace associated to $\{\pm\lambda_p,\pm\overline{\lambda}_p\}$ generates slowly decaying undulations depicted in Figure\,\ref{f:capy}. 

\paragraph{Weakly stable and unstable manifolds} In \cite{PW_14} a center manifold reduction of the spatial dynamics formulation of the stationary, constant coefficient FCH equation was used to classify all solutions that remain close to a bilayer profile as $t\to\pm\infty.$ Excluding the meander modes through a symmetry assumption, a series of normal form transformations were used to recast the four-dimensional pearling center manifold in the form
 \begin{equation}\label{e:PNF}
\begin{aligned}
\dot{C_1}&=\rmi (1+\omega_1\varepsilon) C_1 + C_2+\rmi C_1\big[\alpha_7C_1\bar{C}_1+\alpha_8\rmi(C_1\bar{C}_2-\bar{C}_1C_2)\big],\\
\dot{C_2}&=\rmi (1+\omega_1\varepsilon) C_2 +\rmi C_2\big[\alpha_7C_1\bar{C}_1+\alpha_8\rmi(C_1\bar{C}_2-\bar{C}_1C_2)\big]
+ C_1\left[-\alpha_0\varepsilon+\rmi \alpha_2(C_1\bar{C_2}-\bar{C_1}C_2)\right],
\end{aligned}
\end{equation}
where $C_1$, $C_2\in\C$, the constants $\omega_1, \alpha_j\in\R$, and the conjugate equations are omitted. So that $(C_1, C_2)$ lie on the stable manifold of the bilayer solution we impose the necessary condition  
\[
\lim_{t\to\infty} C_1(t)=\lim_{t\to\infty} C_2(t)=0.
\]
Noting that 
\[
K=\rmi(C_1\bar{C}_2-\bar{C}_1C_2), \qquad H=|C_2|^2-(-\alpha_0\veps+2\alpha_2K)|C_1|^2,
\]
are two first-integrals of \eqref{e:PNF}, we conclude that the stable manifold lies within 
\[
K=H=0.
\]
Introducing the polar coordinates, 
\[
\begin{cases}
C_1=r_1\rme^{\rmi((1+\omega_1\veps)t+\theta_1)},\\
C_2=r_2\rme^{\rmi((1+\omega_1\veps)t+\theta_2)},
\end{cases}
\]
the equation $K=H=0$ can be rewritten as 
\[
\begin{cases}
\theta_1-\theta_2=k\pi, \quad k\in\Z,\\
r_2^2=-\alpha_0\veps r_1^2.
\end{cases}
\]
In addition, we have 
\[
\begin{cases}
2r_1\dot{r}_1=\dot{C}_1\bar{C}_1+C_1\dot{\bar{C}}_1=2(-1)^kr_1r_2,\\
2r_2\dot{r}_2=\dot{C}_2\bar{C}_2+C_2\dot{\bar{C}}_2=2(-1)^k(-\alpha_0\veps)r_1r_2,
\end{cases}
\]
or simply, 
\[
\begin{cases}
\dot{r}_1=(-1)^kr_2,\\
\dot{r}_2=(-1)^k(-\alpha_0\veps)r_1,
\end{cases}
\]
which yields
\[
\begin{pmatrix} r_1(t) \\ r_2(t)\end{pmatrix}=C\begin{pmatrix}1 \\ (-1)^{k+1}\sqrt{-\alpha_0\veps}\end{pmatrix} \rme^{-\sqrt{-\alpha_0\veps}t}.
\]
Since $r_1, r_2\geq0$,  we deduce that $k$ is odd. Plugging the polar coordinates into the first equation of \eqref{e:PNF}, we have
\[
\dot{\theta}_1=\alpha_7r_1^2,
\]
which implies that 
\[
\theta_1=\theta_{1,0}-\frac{\alpha_7C^2}{2\sqrt{-\alpha_0\veps}}\rme^{-2\sqrt{-\alpha \veps }t},
\]
for some $\theta_{1,0}\in\R.$ We conclude that $\lim_{t\to\infty} C_1(t)=\lim_{t\to\infty} C_2(t)=0$ if and only if,
\begin{equation}\label{e:sm}
    \begin{cases}
    C_1=&C_s\exp\left\{-\sqrt{-\alpha_0\veps}t+\rmi\left[(1+\omega_1\veps)t-\frac{\alpha_7C_s^2}{2\sqrt{-\alpha_0\veps}}\exp(-2\sqrt{-\alpha \veps }t)+\theta_s\right] \right\},\\
    C_2=&-\sqrt{-\alpha_0\veps}C_1,
\end{cases}
\end{equation}
where $C_s,\theta_s$ parameterize the stable manifold associated to the bilayer solution.
Similarly, all solutions satisfying $\lim_{t\to-\infty} C_1(t)=\lim_{t\to-\infty} C_2(t)=0$ admit the form
\begin{equation}\label{e:usm}
    \begin{cases}
    C_1=&C_u\exp\left\{\sqrt{-\alpha_0\veps}t+\rmi\left[(1+\omega_1\veps)t+\frac{\alpha_7C_u^2}{2\sqrt{-\alpha_0\veps}}\exp(2\sqrt{-\alpha \veps }t)+\theta_u\right] \right\},\\
    C_2=&\sqrt{-\alpha_0\veps}C_1,
\end{cases}
\end{equation}
where $C_u,\theta_u$ parameterize the unstable manifold associated to the bilayer solution.

It is straightforward to see that the unstable and stable manifolds \eqref{e:sm},\eqref{e:usm} intersects with each other only at the $C_u=C_s=0$, which implies that the only orbit whic is homoclinic orbits to the bilayer is the trivial bilayer solution.
\begin{Lemma}
For $\alpha_0<0$ the system \eqref{e:PNF}, characterizing the stable manifold of the bilayer solution of the constant coefficient stationary equation FCH does not admit any nontrivial orbit that is homoclinic to the bilayer solution.
\end{Lemma}

In the undulating regime, $\alpha_0<0$, the stable and unstable manifolds display the slow, undulated decay associated to defect solutions. However, in the absence of defects there are no stationary homoclinic solutions on the center manifold of the bilayer profile, at least for the system truncated beyond cubic nonlinear terms. In the following section we show that localized perturbations to $\eta_2$ trigger undulated responses that decay on the slow $\veps^{-1/2}$ length scale.  The inhomogeneity in $\eta_2(\tau)$ renders the infinite-dimensional dynamical system non-autonomous, which in turn makes the direct application of center manifold reduction cumbersome, if not impossible. As a remedy we pursue a functional analytic route via a Lyapunov-Schmidt reduction.

\section{Existence of undulated bilayer interfaces}
The center manifold analysis presented in Section \,\ref{s:CM-reduction} highlights the role played by the sign of the pearling parameter $\alpha_0$, defined in \eqref{e:alpha0}. When $\alpha_0>0$ this framework supports the construction of pearled bilayers via a spatial dynamics analysis. In this section, we consider the complementary case, $\alpha_0<0$, and show that the stationary FCH \eqref{e:2sFCH} with inhomogeneous  $\eta_2$ coefficient supports undulated bilayers. More specifically, for a sufficiently strongly localized inhomogeneity we show that the constant width bilayer deforms to a solution of \eqref{e:2sFCH} that has long wavelength width-oscillations that decay like $1/\sqrt{\veps}$ in the fast variables away from the defect. We call these solutions undulated bilayers.
\begin{Definition}
For $\delta>0$ an undulated bilayer with flat interface, denoted $u_n(x;\delta)$, is a solution to the inhomogeneous stationary FCH equation \eqref{e:2sFCH} with $\eta_2$ as in \eqref{eq-eta2-def} subject to the boundary conditions
\begin{equation}
\label{eq:BCs}
\begin{aligned}
\lim_{x_1\to\pm\infty}u_n(x;\delta) &=u_h(x_2/\veps;\veps, \eta_{d,\delta}^\pm),\\
\lim_{x_2\to\pm\infty}u_n(x;\delta) &= u_\infty(\veps,\eta_2(x_1)).
\end{aligned}
\end{equation}
For fixed value of $\eta_d$ the flat bilayer interface $u_h(r;\veps, \eta_d)$ is given by Theorem\,\ref{d:bilayer}, while the far-field value $u_\infty$ is given in \eqref{eq-uinf} and
$\eta_{d,\delta}^\pm:=\lim\limits_{x_1\to\pm\infty}\eta_{d,\delta}(x_1).$
\end{Definition}

As in section 2, we rewrite the inhomogeneous stationary FCH equation \eqref{e:2sFCH} in the inner coordinates 
\beq \label{rescaled}
(t,r)=(\sqrt{\lambda_0}x_1/\veps, x_2/\veps),
\eeq
where it takes the form
\begin{equation}\label{e:in-sFCH-2}
    (\partial_r^2-W^\pprime(u)+\lambda_0\partial_t^2+\veps\eta_1)(\partial_r^2u-W^\prime(u)+\lambda_0\partial_t^2u)+\veps\eta_{d,\delta}(t)W^\prime(u)=\veps\gamma,
\end{equation}
and introduce the modulated bilayer interface  \beq\label{def-mbi}
u_{h,\delta}(t,r;\veps):=u_h(r;\veps,\eta_{d,\delta}(t)).
\eeq
In the inner coordinates $\xi=\xi(t)$ varies on an $O(1)$ length scale, in particular its support is $O(1).$
The modulated bilayer interface satisfies both the $t$-modulated family of ODEs,
\begin{equation}\label{e:MBI}
    \left(\partial_r^2-W^{\prime\prime}(u_{h,\delta})+\veps\eta_1\right)
 \left(\partial_r^2u_{h,\delta}-W^\prime(u_{h,\delta})\right)+\veps\eta_{d,\delta}(t) W^\prime(u_{h,\delta})=\veps\gamma,
\end{equation}
and the boundary conditions \eqref{eq:BCs}, but not the full system \eqref{e:in-sFCH-2}. We construct the undulated bilayer interface as a perturbation of the modulated bilayer interface, 
\begin{equation}\label{e:undulation-sch}
    u(t,r;\delta,\veps)=
    u_{h,\delta}(t,r;\veps)+v(t,r;\delta,\veps).
\end{equation}
Inserting the expansion \eqref{e:undulation-sch} into \eqref{e:in-sFCH-2}
and subtracting the $t$-modulated ODEs from both sides defines the residual
\beq
\label{eq:Residual}
\begin{aligned}
\mrF(v;\delta,\veps):=&\Big(\partial_r^2-W^\pprime(u_{h,\delta}+v)+\lambda_0\partial_t^2+\veps\eta_1\Big)\Big(\partial_r^2v-\big(W^\prime(u_{h,\delta}+v)-W^\prime(u_{h,\delta})\big)+\lambda_0\partial_t^2(u_{h,\delta}+v)\Big)+\\
&\Big(-\big(W^\pprime(u_{h,\delta}+v)-W^\pprime(u_{h,\delta})\big)+\lambda_0\partial_t^2\Big)\Big(\partial_r^2u_{h,\delta}-W^\prime(u_{h,\delta})\Big)+\\
&\veps\eta_{d,\delta}(t)\Big(W^\prime(u_{h,\delta}+v)-W^\prime(u_{h,\delta})\Big),
\end{aligned}
\eeq
which is zero precisely when $u$ satisfies \eqref{e:in-sFCH-2}.

We construct solutions of the system $\mrF(v;\delta,\veps)=0$ through the implicit function theorem. Since the $\delta=0$ problem is homogeneous we have
\beq
\label{eq:IFT1}
\mrF(0,0,\veps)=0.
\eeq
Introducing the linearization of \eqref{e:lODE} about the modulated bilayer, 
\beq
\label{def-cLhd}
\caL_{h,\delta}:=\partial_r^2-W^\pprime(u_{h,\delta}),
\eeq
and the scaled residual of \eqref{e:lODE} at the modulated bilayer 
$$w_{h,\delta}:=\frac{\partial_r^2u_{h,\delta}-W^\prime(u_{h,\delta})}{\veps},$$
it is straightforward to verify that
\beq\label{e:parFv}
 \frac{\partial \mrF}{\partial v}(0,\delta,\veps)= L_{\delta,\veps},
\eeq
where the operator
    \begin{equation}\label{e:Lde}
    \begin{aligned}
     L_{\delta,\veps}= \left(\caL_{h,\delta}+\veps\eta_1+\lambda_0\partial_t^2\right)\left(\caL_{h,\delta}+\lambda_0\partial_t^2\right)+\veps\left(\eta_{d,\delta} W^{\pprime}(u_{h,\delta})-W^\tprime(u_{h,\delta})w_{h,\delta}\right) +
      \lambda_0\partial_t^2u_{h,\delta}W^\tprime(u_{h,\delta}),
    \end{aligned}
    \end{equation}
plays a fundamental role in the analysis.
Our analysis is perturbative from the case $\veps=\delta=0$ for which we have the simple operator studied in section 2,
\beq\label{def-L0}
\frac{\partial \mrF}{\partial v}(0,0,0)=L_0:=(\caL_0+\lambda_0\partial_t^2)^2.
\eeq
Here $\caL_0=\partial_r^2-W^\pprime(u_0)$ is defined in \eqref{e:L0} and  has eigenpairs 
$\{(\lambda_j, \psi_j(r))\}_{j=0}^1$ with $
\lambda_0>\lambda_1=0$ and the remainder of its spectrum strictly negative.

\begin{Remark}
We note that the operator $\caL_{h,\delta}$, defined in \eqref{def-cLhd}, when $\delta=0$, coincides with the operator $\caL_h$ defined in \eqref{e:caLh-def}; that is, $\caL_{h,0}=\caL_h$. 
\end{Remark}

\begin{Remark}
The Lyaponov-Schmidt reduction is markedly simpler in the case $\delta=\veps=0$ as compared to the case $\delta=0,0<\veps\ll1$. The operator $L_0$ admits nontrivial invariant spectral spaces that are separable in $L^2(\R^2)$ as $L_0$. Conversely the linear operator
\[
    \frac{\partial F}{\partial v}(0,0,\veps)=L_{0,\veps}=\Big(\caL_{h,0}+\veps\eta_1+\lambda_0\partial_t^2\Big)\Big(\caL_{h,0}+\lambda_0\partial_t^2\Big)+\veps\Big(\eta_{d,0} W^{\pprime}(u_{h,0})-W^\tprime(u_{h,0})w_{h,0}\Big),
    \]
    does not admit such spaces. More specifically, $L_{0,\veps}$ has the decomposition
    \[
L_{0,\veps}=\mathcal{R}_{h,0}+\lambda_0^2\partial_t^4+\veps\eta_1\lambda_0\partial_t^2+2\lambda_0\partial_t^2\caL_{h,0},
\]
where
\[
\mathcal{R}_{h,0}:=\Big(\caL_{h,0}+\veps\eta_1\Big)\caL_{h,0}+\veps\Big(\eta_{d,0} W^{\pprime}(u_{h,0})-W^\tprime(u_{h,0})w_{h,0}\Big).
\]
The term $2\lambda_0\partial_t^2\caL_{h,0}$ makes the invariant subspaces of $L_{0,\veps}$ non-separable. On the other hand, since $\caL_0=\partial_r^2-W^\pprime(u_0)$ is a Sturm-Louiville operator admitting eigenpairs $\{(\lambda_j, \psi_j(r))\}_{j=0}^\infty$, we may conclude by analytical continuation that $\mathcal{R}_{h,0}$ is a self-adjoint operator admitting eigenpairs $\{(\lambda_{j,\veps}, \psi_{j,\veps}(r))\}_{j=0}^\infty$ with
\[
\begin{aligned}
\lambda_{0,\veps}&=\lambda_0^2+\caO(\veps)>\lambda_{1,\veps}=\caO(\veps),\\ \lambda_{2,\veps}&=\lambda_2^2+\caO(\veps)<\lambda_{3,\veps}=\lambda_3^2+\caO(\veps)<\cdots, \quad \lim_{j\to\infty}\lambda_{j,\veps}=+\infty,
\end{aligned}
\]
and the eigenfunctions $\{\psi_{j,\veps}\}_{j=1}^\infty$ which form a complete orthonormal basis of $L^2(\R)$. 
\end{Remark}

We exploit the fact that for $j=0,1$ the subspaces
$$X_j:=\left\{\psi_j(r)\phi(t)\mid \phi\in L^2(\R)\right\}\subset L^2(\R^2),$$ 
are invariant under the resolvent operator associated to $L_0$.
We introduce the $L_0$-invariant central and hyperbolic subspaces
\[
V_c:=\{\psi_0(r)\phi_0(t)+\psi_1(r)\phi_1(t)\mid \phi_0,\phi_1\in L^2(\R)\}, \quad V_h:=V_c^\perp,
\]
and denote by $P$ the $L^2(\R^2)$ orthogonal projection onto $V_c$ with $Q:=\id-P$. We decompose $v$
\beq\label{v-decomp}
v=v_c+v_h, 
\eeq
where $v_c:=Pv$ and $v_h:=Qv$, and write the residual equation in the projected form
\[
\begin{cases}
P\mrF(v_c+v_h, \delta,\veps)=0,\\
Q\mrF(v_c+v_h, \delta,\veps)=0.\\
\end{cases}
\]
The following Lemma solves the $Q$ equation for $v_h$ given a fixed $v_c$.
\begin{Lemma}
There exists an open neighborhood  $\mathcal{B}_0$ of the origin in $H^4(\R)\times H^4(\R)\times\R\times\R_+$, and a smooth mapping 
\[
H: \mathcal{B}_0\longrightarrow H^4(\R^2)\cap V_h,
\]
such that the decomposition \eqref{v-decomp} of $v$ with $v_c=\psi_0\phi_0+\psi_1\phi_1$ and $v_h=H(\phi_0,\phi_1;\delta,\veps)$ satisfies
\begin{equation}\label{e:H}
    Q\mrF(v_c+H(\phi_0,\phi_1;\delta,\veps);\delta,\veps)=0,
\end{equation}
 for all $(\phi_0,\phi_1,\delta,\veps)\in \mathcal{B}_0.$
\end{Lemma}
\begin{Proof}
We denote $\widetilde{V}_h:=H^4(\R^2)\cap V_h$, take $\delta_0,\veps_0>0$ sufficiently small and introduce the $C^\infty$-smooth mapping
\[
\begin{matrix}
\widetilde{F}:& \widetilde{V}_h\times L^2(\R)\times L^2(\R)\times (0,\delta_0)
\times(0,\veps_0) &\longrightarrow   & V_h \\
&(v_h,\phi_0,\phi_1;\delta,\veps) &\longmapsto & QF(\psi_0\phi_0+\psi_1\phi_0+v_h;\delta,\veps).
\end{matrix}
\]
From its construction $\widetilde{F}(0,0,0;0,0)=0.$ It is easy to verify that
$$ \frac{\partial \widetilde{F}}{\partial v_h}(0,0,0;0,0)=QL_0Q: \widetilde{V_h}\to V_h.$$
Moreover, since $Q$ is a spectral projection for $L_0$, $\caL_0$ is strictly negative on $V_h$, and $\lambda_0>0$ we deduce that $\sigma(QL_0Q)$ is bounded away from zero and $QL_0Q$ has a bounded inverse. We may apply the implicit function theorem to solve $\widetilde{F}=0,$ and
concluded the proof.
\end{Proof}

With this reduction of $v$ the $P$ equation may be written in the form
\begin{equation}\label{e:LS}
\begin{cases}
\left\langle \psi_0, \mrF\bigl(\psi_0\phi_0+\psi_1\phi_1+H(\phi_0,\phi_1;\delta,\veps);\delta,\veps\bigr)\right\rangle_{L^2(\R)}=0, \\
\langle \psi_1,\mrF(\psi_0\phi_0+\psi_1\phi_1+H(\phi_0,\phi_1;\delta,\veps);\delta,\veps)\rangle_{L^2(\R)}=0,
\end{cases}
\end{equation}
where the left hand sides depend on $t$ via $\phi_0$, $\phi_1$ and $H$. 
 The inhomogeneity in $\eta_2$ does not break the $x_2$ even parity of the stationary FCH equation. Without loss of generality we restrict ourselves to functions with even parity in $r$ for each fixed $t$. More specifically, we introduce
\[
v\in L^2_{\mathrm{even}}(\R^2):=\{v\in L^2(\R^2)\mid v(r,t)=v(-r,t)\},
\]
and remark that $\phi_1(t)\equiv 0$ for $v\in L^2_{\mathrm{even}}(\R^2)$ since $\psi_1$ has odd parity. The second equation in \eqref{e:LS} holds trivially since $F(v;\delta,\veps)$ has even parity, and the first equation simplifies to
\begin{equation}\label{e:K}
K(\phi_0;\delta,\veps):= \left\langle \psi_0, F\bigl(\psi_0\phi_0+H(\phi_0,0;\delta,\veps);\delta,\veps\bigr)\right\rangle_{L^2(\R)}=0.
\end{equation}
In the sequel we fix $\phi_1=0$ and drop references to it. We
 assume that $|\delta|, \veps>0$ are sufficiently small and use the contraction mapping principle to construct the solution $\phi_0$ of \eqref{e:K} and identify is leading order form, establishing Theorem\,\ref{thm-main}.

The map $K$ is smooth, $$K:\widetilde{\mathcal{B}_0}:=\{(\phi_0;\delta,\veps)\in H^4(\R)\times \R\times[0,\infty)\mid (\phi_0, 0;\delta,\veps)\in\mathcal{B}_0\}\mapsto L^2(\R),$$ and admits the expansion
\begin{equation}\label{e:Kexp}
    K(\phi_0;\delta,\veps)=K(0;\delta,\veps)+\frac{\partial K}{\partial\phi_0}(0;\delta,\veps)\phi_0+\caO(\|\phi_0\|_{H^4(\R)}^2).
\end{equation}

The equality \eqref{eq:IFT1} suggests that the leading order term in \eqref{e:Kexp} is small, correspondingly we introduce
\begin{equation}\label{def-K0}
  K_0(\cdot;\delta,\veps):=\frac{K(0;\delta,\veps)}{\delta\veps}=\frac{\langle\psi_0,\mrF(H(0;\delta,\veps);\delta,\veps)\rangle_{L^2(\R)}}{\delta\veps}\in L^2(\R).
\end{equation}

\begin{Lemma}
\label{l:K0}
The function $K_0\in L^2(\R)$ defined in \eqref{def-K0} has the leading order expansion
\beq
\label{K0-exp}
K_0(t;\delta, \veps)=
\Big\langle \psi_0,W^\prime(u_0) \Big\rangle_{L^2(\R)}\Big( \xi^{(4)}(t)+2 \xi^{\pprime}(t)\Big) +\caO\Big(|\veps|+|\delta\veps|\Big),
\eeq
where the inhomogeniety $\xi$ is defined in \eqref{def-xi} and the error is measured in $L^2(\R).$
\end{Lemma}
\begin{Proof}
It is convenient to write the modulated bilayer interface $u_{h,\delta}$ of \eqref{def-mbi} as a perturbation of the flat bilayer; that is, the Taylor expansion of $u_{h,\delta}$ with respect to $\delta$, which reads
\[
u_{h,\delta}(t,r)=u_{h,0}(r)+\delta\veps \xi(t)\Big(L(\veps,0)^{-1}W^\prime(u_0)\Big)(r)+\caO(|\delta\veps|^2),
\]
where the one-dimensional operator $L(\veps,\lambda)$ is defined in \eqref{e:1eigen}. Similarly, the residual $\mrF$, from \eqref{eq:Residual}, admits the expansion
\begin{equation}\label{e:Fexp}
    F(v;\delta,\veps)=F(0,\delta,\veps)+L_{\delta,\veps}v+\caO(\|v\|^2_{H^4(\R)}),
\end{equation}
where the modulated residual is given by
\begin{equation}\label{e:F0exp}
    \begin{aligned}
F(0,\delta,\veps)=&\left(\partial_r^2-W^\pprime(u_{h,\delta})+\lambda_0\partial_t^2+\veps\eta_1 \right)\left(\lambda_0\partial_t^2u_{h,\delta} \right)+\lambda_0\partial_t^2\left(\partial_r^2u_{h,\delta}-W^\prime(u_{h,\delta}) \right)\\
=& \veps\delta F_{0,\veps} +\caO(|\delta\veps|^2).
\end{aligned}
\eeq
Here, using the expansion of $u_{h,\delta}$, we have introduced
\beq\label{def-Foe}
F_{0,\veps}:=\lambda_0^2\xi^{(4)}(t)\Big(L^{-1}(\veps,0)W^\prime(u_0) \Big)(r)+\lambda_0\xi^{\pprime}(t)\Big[(2\caL_{h,0}+\veps\eta_1)L^{-1}(\veps,0)W^\prime(u_0)\Big](r).
\end{equation}
The result \eqref{e:H} holds with $\phi_0=0$, which implies that 
\begin{equation}\label{e:H0a}
    Q\mrF\bigl(H(0;\delta,\veps);\delta,\veps\bigr)=0.
\end{equation}
Using \eqref{e:Fexp} and \eqref{e:F0exp} to expand \eqref{e:H0a},  
yields the asymptotic result
\begin{equation}\label{e:Hexp}
    H(0;\delta,\veps)=-\delta\veps H_0(\veps) + \caO\Big(|\delta\veps|^2\Big),
\end{equation}
where we have introduced
\beq\label{def-H0}
H_0(\veps):=\Big(QL_{0,\veps}Q\Big)^{-1}\Big(QF_{0,\veps}\Big).
\eeq
Combining the expansions for $F_{\delta,\veps}$ and $H$ (\ref{e:Fexp}-\ref{e:Hexp}) we simplify $K_0$ as
\begin{equation}\label{e:fexp}
\begin{aligned}
K_0=
 &\Big\langle \psi_0,F_{0,\veps}-L_{0,\veps}H_0\Big\rangle_{L^2(\R)}+\caO\Big(|\delta\veps|\Big)\\
=&\Big\langle \psi_0,F_{0,0}-L_{0,0}\Big(QL_{0,0}Q\Big)^{-1}\Big(QF_{0,0}\Big)\Big\rangle_{L^2(\R)}+\caO\Big(|\veps|+|\delta\veps|\Big).
\end{aligned}
\eeq
Since $Q$ commutes with $L_{0,0}$ and $Q\psi_0=0$ the second term in the inner product is zero and we are left with the inner product with $F_{0,0}.$ From \eqref{def-Foe} with $\veps=0$ we arrive at \eqref{def-K0}. Since the derivatives of $\xi$ are compactly supported $K_0$ is too, at least at leading order. 
\end{Proof}

To simplify the second term on the right-hand side of \eqref{e:Kexp} we introduce the constant coefficient operator
\beq\label{cG-def}
\cG:=(1+\partial_t^2)^2+\veps c_1(1+\partial_t^2)-4\veps\alpha_0,
\eeq
that acts purely through the tangential variable. Here the constant $c_1\in\R$ is defined by
$$c_1:=\frac{\eta_1-2\langle \psi_0, W^\tprime(u_0)u_1\psi_0 \rangle_{L^2(\R)}}{\lambda_0},$$ 
where $\psi_0$ is normalized to have $L^2$ norm 1 and $u_1$ is the $\veps$-scaled leading order term of $u_h-u_0$ as defined in \eqref{e:u1w0}. Significantly the quantity $\alpha_0$ is the key bifurcation parameter defined in \eqref{e:alpha0} that establishes the positivity of $\cG$ for $\alpha_0<0.$

\begin{Lemma}
\label{l:pK0}
\beq \label{e:pK_0}
\frac{\partial K}{\partial \phi_0}(0;\delta,\veps)\phi_0 = \lambda_0^2 \cG \phi_0 + \caO\Big((\delta\veps+\veps^2)\|\phi_0\|_{H^4(\R)}\Big).
\eeq
\end{Lemma}
\begin{Proof}
Recalling (\ref{e:parFv}, \ref{e:K}), a straightforward calculation shows that
\begin{equation}\label{e:K1exp}
\begin{aligned}
\frac{\partial K}{\partial \phi_0}(0;\delta,\veps)\phi_0=&\left\langle \psi_0, L_{\delta,\veps}\left(\psi_0\phi_0+\frac{\partial H}{\partial \phi_0}(0;\delta,\veps)\phi_0 \right)\right\rangle_{L^2(\R)}\\
=&\left\langle \psi_0, L_{0,\veps}\left(\psi_0\phi_0+\frac{\partial H}{\partial \phi_0}(0;0,\veps)\phi_0 \right)\right\rangle_{L^2(\R)}+\caO\bigg(|\delta\veps|\|\phi_0\|_{H^4(\R)}\bigg),
\end{aligned}
\end{equation}
where the latter equality is simply the leading order expansion in $\delta$.
From \eqref{e:Lde} with $\delta=0$, the operator $L_{0,\veps}$ takes the form
\begin{equation}\label{e:Lveps}
   \begin{aligned}
   L_{0,\veps}v=&\Big(\caL_{h,0}+\veps\eta_1+\lambda_0\partial_t^2\Big)\Big(\caL_{h,0}+\lambda_0\partial_t^2\Big)v+\veps\Big(\eta_{d,0} W^{\pprime}(u_{h,0})-W^\tprime(u_{h,0})w_{h,0}\Big)v\\
   =&L_0 v+\veps L_1v +\caO(\veps^2),
\end{aligned} 
\end{equation}
where the last expression gives the Taylor expansion of $L_{0,\veps}$ in $\veps$ with  $L_0=(\caL_0+\lambda_0\partial_t^2)^2$ as in \eqref{def-L0} and the first order operator
\[
L_1v:=\Big(\eta_1-W^\tprime(u_0)u_1\Big)\Big(\caL_0+\lambda_0\partial_t^2\Big)v-\Big(\caL_0+\lambda_0\partial_t^2\Big)\Big( W^\tprime(u_0)u_1v\Big)+\Big(\eta_{d,0} W^{\pprime}(u_0)-W^\tprime(u_0)w_0\Big)v.
\]
Using \eqref{cG-def} and the expansion \eqref{e:Lveps}, the first term on the right-hand side of \eqref{e:K1exp} is expressed as
\begin{equation}\label{e:K1exp-1}
\begin{aligned}
\Big\langle \psi_0, L_{0,\veps}\left(\psi_0\phi_0 \right)\Big\rangle_{L^2(\R)}=&
\Big\langle \psi_0, \Big[L_0+\veps L_1\Big]\left(\psi_0\phi_0 \right)\Big\rangle_{L^2(\R)}+\caO\Big(\veps^2\|\phi_0\|_{H^4(\R)}\Big),
\\
=&
\Big\langle \psi_0, \left[\lambda_0^2\Big(1+\partial_t^2\Big)^2+\lambda_0\veps\Big(\eta_1-2W^\tprime(u_0)u_1\Big)\Big(1+\partial_t^2\Big)+
\right.\\
&\left.\veps\Big(\eta_{d,0} W^{\pprime}(u_0)-W^\tprime(u_0)w_0\Big)\right]\left(\psi_0\phi_0 \right)\Big\rangle_{L^2(\R)}+\caO\Big(\veps^2\|\phi_0\|_{H^4(\R)}\Big),
\\
=&\lambda_0^2\,\mathcal{G}\phi_0 +\caO\Big(\veps^2\|\phi_0\|_{H^4(\R)}\Big).
\end{aligned}
\end{equation}

To estimate the second, lower-order term on the right-hand side of \eqref{e:K1exp}  we first expand $\frac{\partial H}{\partial \phi_0}(0; 0,\veps)$.  When $\phi_1=0$ and $\delta=0$ equation \eqref{e:H} reduces to 
\[
QF(\psi_0\phi_0+H(\phi_0; 0,\veps),0,\veps)=0.
\]
Linearizing this relation with respect to $\phi_0$ at $\phi_0=0$, yields
\begin{equation}\label{e:H_phi-0}
    QL_{0,\veps}\left(\psi_0\phi_0+\frac{\partial H}{\partial \phi_0}(0; 0,\veps)\phi_0 \right)=0.
\end{equation}
Since $Q$ commutes with $L_0$ and $Q\psi_0\phi_0=0$, applying $Q$ to the expansion \eqref{e:Lveps} yields 
$$
QL_{0,\veps}(\psi_0\phi_0)=\veps QL_1(\psi_0\phi_0)+\caO(\veps^2\|\phi_0\|_{H^4(\R)}).$$
We expand $\partial_{\phi_0} H$ as 
\begin{equation}\label{e:Hsch}
    \frac{\partial H}{\partial \phi_0}(0; 0,\veps)=H_0+\veps H_1+\caO(\veps^2),
\end{equation}
and plug this into \eqref{e:K1exp}, and equate orders of $\veps$. At first and second order we find
\begin{equation}\label{e:Hdeco}
    \begin{cases}
    QL_0H_0\phi_0=0,\\
    QL_0H_1\phi_0+QL_1H_0\phi_0+QL_1(\psi_0\phi_0)=0,
    \end{cases}
\end{equation}
for any $\phi_0\in H^4(\R)$.
The operator $\frac{\partial H}{\partial \phi_0}(0; 0,\veps)$ maps $H^4(\R)$ to $V_h\cap H^4(\R^2)$ while $QL_0Q$ is $\veps$-uniformly invertible, on $V_h$ which is the range of $Q$. 
Since the expansion holds for all $\phi_0\in H^4$ the first equation in \eqref{e:Hdeco} implies
\begin{equation}\label{e:H0b}
    H_0 =0,
\end{equation}
and the second equation in \eqref{e:Hdeco} reduces to
\[
QL_0Q H_1\phi_0=-QL_1(\psi_0\phi_0),
\]
or, equivalently,
\begin{equation}\label{e:H1}
    H_1\phi_0=-\Big(QL_0Q\Big)^{-1}\Big( QL_1(\psi_0\phi_0)\Big).
\end{equation}
With these formulations the expansion  \eqref{e:Hsch}, reduces to
\begin{equation}\label{e:Hsch-ex}
    \frac{\partial H}{\partial \phi_0}(0; 0,\veps)\phi_0=-\veps \Big(QL_0Q\Big)^{-1}\Big( QL_1(\psi_0\phi_0)\Big)+\caO(\veps^2\|\phi_0\|_{H^4(\R)})\in V_h,
\end{equation}
and together with \eqref{e:Lveps}
we rewrite the second term on the right-hand side of \eqref{e:K1exp} as
\beq\label{e:pH0-est}
\begin{aligned}
\left\langle \psi_0,  L_{0,\veps}\left(\frac{\partial H}{\partial \phi_0}(0,;0,\veps)\phi_0 \right)\right\rangle_{L^2(\R)}=
&-\veps\left\langle \psi_0, L_0\Big(QL_0Q\Big)^{-1}\Big( QL_1(\psi_0\phi_0)\Big)\right\rangle_{L^2(\R)}+\caO(\veps^2\|\phi_0\|_{H^4(\R)}),\\
=&-\veps\left\langle \psi_0,  QL_1\Big(\psi_0\phi_0\Big)\right\rangle_{L^2(\R)}+\caO(\veps^2\|\phi_0\|_{H^4(\R)}), \\
= & \caO\bigg(\veps^2\|\phi_0\|_{H^4(\R)}\bigg).
\end{aligned}
\eeq
Together \eqref{e:K1exp-1} and \eqref{e:pH0-est}  establish \eqref{e:pK_0}.
\end{Proof}

\begin{Remark}
\label{r:mbbL-L}
In section 2 the stationary FCH equation was rewritten as an infinite-dimensional dynamical system \eqref{e:idds},
\[
\dot{V}=\bbL(\veps) V +\bbF(V,\veps).
\]
The eigenvalue problem $\bbL(\veps) V=\lambda V$ is equivalent to \eqref{e:1eigen},
\[
L(\veps,\lambda)v:=\bigg[ \left(\caL_h+\veps\eta_1+\lambda_0\lambda^2\right)\left(\caL_h+\lambda_0\lambda^2\right)+\veps(\eta_d W^{\pprime}(u_h)-W^\tprime(u_h)w_h)\bigg]v=0,
\]
which, restricted to the central mode $v=\psi_0(r)$, leads to a quartic polynomial in $\lambda$,
\[
p(\lambda):=\langle\psi_0, L(\veps,\lambda)\psi_0\rangle_{L^2(\R)}.
\]
The characteristic polynomial of $\cG$ equals $\lambda_0^{-2}p(\lambda)$ to $O(\veps)$. 
\end{Remark}

With Lemmas\,\ref{l:K0} and \ref{l:pK0} we may rewrite the expansion \eqref{e:Kexp} of $K$ as
\beq\label{Kexp2}
K(\phi_0;\delta,\veps) = \delta\veps K_0(\cdot;\delta,\veps)+ \lambda_0^2 \cG \phi_0 + \cR(\phi_0;\delta,\veps),
\eeq
where $K_0\in L^2(\R)$, with expansion given in \eqref{K0-exp}, is independent of $\phi_0$ and the remainder satisfies
\beq
\label{e:R-est}
\cR(\phi_0;\delta,\veps)=\caO\left((\veps^2+|\delta|\veps)\|\phi_0\|_{H^4}
+\|\phi_0\|_{H^4}^2\right).
\eeq
In particular $\phi_0$ satisfies $K(\phi_0;\delta,\veps)=0$ if and only if
\begin{equation}\label{e:LS-2}
    \mathcal{G}\phi_0=-\delta\veps\lambda_0^{-2}K_0(\cdot;\delta,\veps)+\cR(\phi_0;\delta,\veps),
\end{equation}
where we have rescaled the remainder. 

The operator $\cG$, given in \eqref{cG-def} is constant coefficient. Its Green's function $G$ can be determined explicitly, see for example \cite{choksi-21}.  Indeed, $G\in H^3(\R)$ has two continuous derivatives and satisfies
\begin{equation}\label{e:G}
    G(t):=\frac{\rme^{-B|t|}\Big[A\cos(At)+B(\chi(t)-\chi(-t))\sin(At) \Big]}{4AB(A^2+B^2)}=\frac{AE_e(Bt)\cos(At)+BE_o(Bt)\sin(At) }{4AB(A^2+B^2)},
\end{equation}
where $\chi$ is the standard step function and
\beq\label{e:AB}
\begin{aligned}
A=&\sqrt{\frac{\sqrt{1+\veps c_1-4\veps\alpha_0}}{2}+\frac{2+\veps c_1}{4}}=1+\caO(\veps), \\
B=&\sqrt{\frac{\sqrt{1+\veps c_1-4\veps\alpha_0}}{2}-\frac{2+\veps c_1}{4}}=\sqrt{-\alpha_0\veps}+\caO(\veps^{3/2}),\\
E_e(x)=&\begin{cases}
e^{-x}, &x\geq0,\\
e^{x}, &x<0,
\end{cases} 
\qquad
E_o(x)=\begin{cases}
e^{-x}, &x>0,\\
-e^{x}, &x<0.
\end{cases}
\end{aligned}
\eeq

Inverting $\cG$, the relation \eqref{e:LS-2} reduces to a fixed point problem 
\[
\phi_0=T(\phi_0),
\]
where the map
\begin{equation}
    \begin{matrix}
    T: & H^4(\R) & \longrightarrow & H^4(\R)\\
    & \phi_0 & \longmapsto & G*(-\delta\veps\lambda_0^{-2}K_0+\cR),
    \end{matrix}
\end{equation}
is defined through convolution with the Green's function $G$ over $t\in\R$.

\begin{Lemma}\label{l:G*f}
Fix $\veps_0, \delta_0>0$ sufficiently small, then for each $q\in[1,2]$ there exists $C_0>0$ such that for all $\veps\in(0,\veps_0)$ and $\delta\in(-\delta_0,\delta_0)$ we have the estimate
\beq
\label{e:G-est}
\|G*f\|_{H^4(\R)}\leq C_0\veps^{-\frac{5q-2}{4q}}\left(\|f\|_{L^q(\R)}+\|f\|_{L^2(\R)}\right),
\eeq
for all $f\in L^q(\R)\cap L^2(\R)$.
\end{Lemma}
\begin{Proof}
From Young's convolution inequality we may estimate
\[
\|G*f\|_{L^2(\R)}\leq 
\|G\|_{L^p(\R)}\|f\|_{L^q(\R)},
\]
for conjugate exponents $p^{-1}+q^{-1}=\frac32.$ In particular we need $q\in[1,2]$ so that $p\geq 1.$ 
From the expression \eqref{e:G} and asymptotic expansions \eqref{e:AB} we have
\[\|G\|_{L^p(\R)}\leq \frac{A+B}{4AB(A^2+B^2)}\left\|\rme^{-B|\cdot|}\right\|_{L^p(\R)}\leq\frac{A+B}{2AB(A^2+B^2)}  B^{-p^{-1}} \leq C_0 \veps^{-\frac{5q-2}{4q}}.\]
Since $G\in H^3(\R)$ we may take two derivatives of $G$ point-wise and estimate
$$ \|\partial_t^2 G \|_{L^p} \leq C_0 \veps^{-\frac{5q-2}{4q}}.$$
This estimate, and the convolution identity $\partial_t^2( G*f) =(\partial_t^2 G)*f$, implies an $H^2(\R)$ bound on $G*f$ of the form \eqref{e:G-est}. Expanding the operator $\cG$ in \eqref{cG-def}, denoting $u=G*f$, and taking the $L^2$ norm of both sides we have
$$\|\partial_t^4 u\|_{L^2} \leq C \left(\|\partial_t^2u\|_{L^2} +\|u\|_{L^2} + \|f\|_{L^2}\right).$$
Here the constant $C$ may be chosen independent of $\veps$ and $\delta$ sufficiently small. This extends the estimate to $H^4(\R)$ as in \eqref{e:G-est}.
\end{Proof}

\begin{Lemma}\label{l:GK0}
Let the function $K_0$ be as defined in \eqref{def-K0}. Fix $\veps_0, \delta_0$, sufficiently small, then there exists a constant $C_1>0$ such that
\beq\label{e:oK0-est}
\|G*K_0\|_{H^4(\R)} \leq C_1 \veps ^{-3/4},
\eeq
holds for all $(\delta,\veps)\in(-\delta_0,\delta_0)\times(0,\veps_0).$
\end{Lemma}
\begin{Proof}
From the expansion \eqref{K0-exp} the leading order behaviour of $K_0(\cdot; \delta,\veps)$, which we denote by $\oK_0$, is given by
$$
\oK_0:=\Big\langle \psi_0,W^\prime(u_0) \Big\rangle_{L^2(\R)}\Big( \xi^{(4)}(t)+2 \xi^{\pprime}(t)\Big),$$ 
while the remainder $\oK_1:=K_0-\oK_0$ satisfies
$$\oK_1=\caO(\veps+\delta\veps),$$
in $L^2(\R)$.
The function $\oK_0$ inherits compact support from the perturbation $\xi$, however the estimate \eqref{e:oK0-est} only requires $\oK_0\in L^1(\R)\cap L^2(\R)$, for which Lemma \ref{l:G*f} with $q=1$ implies
\begin{equation}\label{e:oK0-0}
    \|G*\oK_0\|_{H^4(\R)} \leq C_0\|\oK_0\|_{L^1(\R)} \veps ^{-3/4}.
\end{equation}
Similarly, Lemma \ref{l:G*f} with $q=2$ applied to $G*\oK_1$, yields
\begin{equation}\label{e:oK0-1}
    \|G*\oK_1\|_{H^4(\R)} \leq \widetilde{C}_0 (1+|\delta|),
\end{equation}
for some $\widetilde{C}_0>0$. We take $C_1:=\max\{C_0\|\oK_0\|_{L^1(\R)},\widetilde{C}_0 (1+|\delta_0|)\veps_0^{3/4} \}$ and  deduce \eqref{e:oK0-est}. 
\end{Proof}
In the sequel we fix 
$$q>3/4$$ 
and rescale $\delta$ and $\phi_0$, introducing
\begin{equation}\label{e:scaling}
    \delta=\widetilde{\delta}\veps^q, \quad \phi_0=\delta\veps^{1/4}\widetilde{\phi}_0,
\end{equation}
and the associated operators
\begin{equation}\label{e:tilde-TM}
     T(\phi_0)=\delta\veps^{1/4}\widetilde{T}(\widetilde{\phi}_0), \quad 
     \cR(\phi_0;\delta,\veps)=\delta\veps^{1/4}\widetilde{\cR}(\widetilde{\phi}_0;\widetilde{\delta},\veps). 
\end{equation}
The rescaled fixed point equation takes the form
\begin{equation}\label{e:tildeT}
    \widetilde{\phi}_0=\widetilde{T}(\widetilde{\phi}_0)=G*\Big(-\veps^{3/4}\lambda_0^{-2}K_0+\widetilde{\cR} \Big),
\end{equation}
where, expanding \eqref{e:R-est}, the residual term satisfies 
\begin{equation}\label{e:tildeM}
    \|\widetilde{\cR}(\widetilde{\phi}_0;\widetilde{\delta},\veps)\|_{L^2(\R)}=\caO(|\veps|^2\|\widetilde{\phi}_0\|_{H^4}+|\widetilde{\delta}\veps^{q+1}|\,\|\widetilde{\phi}_0\|_{H^4}+|\widetilde{\delta}\veps^{q+1/4}|\,\|\widetilde{\phi}_0\|_{H^4}^2)
\end{equation}
The following proposition establishes that $\widetilde{T}$ is a contraction mapping.
\begin{Proposition}
\label{p:contraction}
There exist $\delta_0,\veps_0, R>0$ such that for any given $(\widetilde{\delta},\veps)\in(-\delta_0,\delta_0)\times(0,\veps_0)$, the mapping 
\[
\widetilde{T}: B_0(R):=\{\phi\in H^4(\R)| \|\phi\|_{H^4(\R)}\leq R\}\longrightarrow B_0(R)
\]
is a well-defined contraction and admits a unique fixed point, $\widetilde{\phi}_0^*(t;\widetilde{\delta},\veps)$.
\end{Proposition}

\begin{Proof}
From \eqref{e:oK0-est} there exist $\delta_1,\veps_1, R_1>0$ such that, for any $|\widetilde{\delta}|<\delta_1, 0<\veps<\veps_1$,
\begin{equation}\label{e:Gf}
    \|G*(\veps^{3/4}\lambda_0^{-2}K_0)\|_{H^4(\R)}\leq R_1.
\end{equation}
Similarly, from Lemma\,\ref{l:G*f} and the estimate \eqref{e:tildeM}, we conclude that 
there exist $\delta_2,\veps_2>0, C_2\geq1$ such that, for any $|\widetilde{\delta}|<\delta_2, 0<\veps<\veps_2$,
\begin{equation}\label{e:GM}
    \|G*\widetilde{\cR}(\widetilde{\phi}_0)\|_{H^4(\R)}\leq C_2\Big( |\veps|\|\widetilde{\phi}_0\|_{H^4}+|\widetilde{\delta}\veps^{q}|\|\widetilde{\phi}_0\|_{H^4}+|\widetilde{\delta}\veps^{q-3/4}|\|\widetilde{\phi}_0\|_{H^4}^2\Big).
\end{equation}
We choose $\delta_3,\veps_3>0$ sufficiently small that
\begin{equation}\label{e:radi}
    C_2\max\limits_{|\widetilde{\delta}|<\delta_3, 0<\veps<\veps_3}\{\veps, |\widetilde{\delta}\veps^q|, 4R_1|\widetilde{\delta}\veps^{q-3/4}| \} \leq\frac{1}{4}.
\end{equation}
Combining (\ref{e:Gf})-(\ref{e:radi}) and defining
\[
\delta_4:=\min\{\delta_1,\delta_2,\delta_3\},\quad \veps_4:=\min\{\veps_1,\veps_2,\veps_3\},\quad
R=4R_1,
\]
Then for any $|\widetilde{\delta}|<\delta_4$, $0<\veps<\veps_4$, and all
$\widetilde{\phi}_0\in B_0(R)$ we have the bound
\[
\begin{aligned}
\|\widetilde{T}(\widetilde{\phi}_0)\|_{H^4(\R)}\leq & \|G*(\veps^{3/4}\lambda_0^{-2}\oK_0)\|_{H^4(\R)}+ \|G*\widetilde{\cR}(\widetilde{\phi}_0)\|_{H^4(\R)}\\
\leq & R_1+\frac{1}{4}\Big(2\|\widetilde{\phi}_0\|_{H^4}+\frac{1}{4R_1}\|\widetilde{\phi}_0\|_{H^4}^2\Big)\\
\leq &R_1+\frac{1}{4}\Big(2R+\frac{R^2}{4R_1}\Big)=R.
\end{aligned}
\]
This establishes that $\widetilde{T}: B_0(R)\to B_0(R)$. Taking $\widetilde{\phi}_0, \widetilde{\varphi}_0\in B_0(R)$, we bound 
$\|\widetilde{T}(\widetilde{\phi}_0)-\widetilde{T}(\widetilde{\varphi}_0)\|_{H^4(\R)}$ in terms of $\|\widetilde{\phi}_0-\widetilde{\varphi}_0\|_{H^4(\R)}$. Since $K_0$ is independent of $\phi_0$ we have
\begin{equation}\label{e:con-0}
    \widetilde{T}(\widetilde{\phi}_0)-\widetilde{T}(\widetilde{\varphi}_0)=
G*\Big(\widetilde{\cR}(\widetilde{\phi}_0)-\widetilde{\cR}(\widetilde{\varphi}_0)\Big).
\end{equation}
From the expansion \eqref{Kexp2} 
and the rescaling in (\ref{e:tilde-TM}), we conclude that 
\begin{equation}
    \cR=\overbrace{-\mathcal{G}\phi_0+\frac{1}{\lambda_0^2}\frac{\partial K}{\partial \phi_0}(0;\delta,\veps)\phi_0}^{\cR_1}+\overbrace{\frac{1}{\lambda_0^2}\Big(K(\phi_0;\delta,\veps)-K(0;\delta,\veps)-\frac{\partial K}{\partial \phi_0}(0,\delta,\veps)\phi_0\Big)}^{\cR_2}
\end{equation}
where $\cR_1$ is linear in $\phi_0$ while $\cR_2$ is genuinely nonlinear in $\phi_0$. We rescale $\widetilde{\cR}_{1\backslash2}$ as in \eqref{e:tilde-TM}. The linearity of $\cR_1$, and the estimates \eqref{e:GM} and \eqref{e:radi} show that 
\begin{equation}\label{e:con-1}
    \|G*\Big(\widetilde{\cR}_1(\widetilde{\phi}_0)-\widetilde{\cR}_1(\widetilde{\varphi}_0)\Big)\|_{H^4(\R)}\leq \frac{1}{2}\|\widetilde{\phi}_0-\widetilde{\varphi}_0\|_{H^4(\R)}.
\end{equation}
For the rescaled nonlinear term $\cR_2$, we claim that there exists $C_3>0$ such that 
\begin{equation}\label{e:con-2}
   \left \|G*\Big(\widetilde{\cR}_2(\widetilde{\phi}_0)-\widetilde{\cR}_2(\widetilde{\varphi}_0)\Big)\right\|_{H^4(\R)}\leq C_3\widetilde{\delta}\veps^{q-3/4}\|\widetilde{\phi}_0-\widetilde{\varphi}_0\|_{H^4(\R)}.
\end{equation}
Exploiting the integral remainder form of Taylor's expansion  we write 
\[
\cR_2(\phi_0)=\frac{1}{\lambda_0^2}\int_0^1\Big(\frac{\partial K}{\partial \phi_0}(t\phi_0;\delta,\veps)-\frac{\partial K}{\partial \phi_0}(0;\delta,\veps)\Big)\phi_0\rmd t,
\]
which for the rescaled operator leads to
\begin{equation}\label{e:Rsub}
\begin{aligned}
\widetilde{\cR}_2(\widetilde{\phi}_0)-\widetilde{\cR}_2(\widetilde{\varphi}_0)=&\frac{1}{\lambda_0^2}\left[\int_0^1\Big(\frac{\partial K}{\partial \phi_0}(t\delta\veps^{1/4}\widetilde{\phi}_0;\delta,\veps)-\frac{\partial K}{\partial \phi_0}(0;\delta,\veps)\Big)\Big(\widetilde{\phi}_0-\widetilde{\varphi}_0\Big)\rmd t+\right.\\
&\left.\int_0^1\Big(\frac{\partial K}{\partial \phi_0}(t\delta\veps^{1/4}\widetilde{\phi}_0;\delta,\veps)-\frac{\partial K}{\partial \phi_0}(t\delta\veps^{1/4}\widetilde{\varphi}_0;\delta,\veps)\Big)\widetilde{\varphi}_0\rmd t\right],\\
=&\frac{\delta\veps^{1/4}}{\lambda_0^2}\left[\int_0^1 t\Big(\int_0^1\frac{\partial^2 K}{\partial \phi_0^2}(ts\delta\veps^{1/4}\widetilde{\phi}_0;\delta,\veps)\widetilde{\phi}_0\rmd s\Big)\Big(\widetilde{\phi}_0-\widetilde{\varphi}_0\Big)\rmd t+\right.\\
&\left.\int_0^1t\Big(\int_0^1\frac{\partial^2 K}{\partial \phi_0^2}(t\delta\veps^{1/4}(\widetilde{\varphi}_0+s(\widetilde{\phi}_0-\widetilde{\varphi_0}));\delta,\veps)\Big(\widetilde{\phi}_0-\widetilde{\varphi}_0\Big)\rmd s\Big)\widetilde{\varphi}_0\rmd t\right].
\end{aligned}
\end{equation}
From \eqref{e:Rsub} we derive the existence of a constant $\widetilde{C}_3$ such that
\[
\|\widetilde{\cR}_2(\widetilde{\phi}_0)-\widetilde{\cR}_2(\widetilde{\varphi}_0)\|_{L^2(\R)}\leq \widetilde{C}_3\delta\veps^{1/4}\|\widetilde{\phi}_0-\widetilde{\varphi}_0\|_{H^4(\R)},
\]
which, together with the Young's inequality in Lemma \ref{l:G*f} with $q=2$, leads to the inequality \eqref{e:con-2}.

Setting
\[
\veps_0:=\veps_4,\qquad \delta_0:=\min\left\{\delta_4, \frac{\veps_0^{3/4-q}}{4C_3}\right\},
\]
the estimates (\ref{e:con-0}-\ref{e:con-2}) imply that for any $|\widetilde{\delta}|<\delta_0, 0<\veps<\veps_0$, $\widetilde{T}$ is a contraction in the sense that 
\[
\|\widetilde{T}(\widetilde{\phi}_0)-\widetilde{T}(\widetilde{\varphi}_0)\|_{H^4(\R)}=
\|G*\Big(\widetilde{\cR}(\widetilde{\phi}_0)-\widetilde{\cR}(\widetilde{\varphi}_0)\Big)\|_{H^4(\R)}
\leq  \frac{3}{4}\|\widetilde{\phi}_0-\widetilde{\varphi}_0\|_{H^4(\R)}.
\]
Since $\widetilde T$ is a strict contraction from $B_0(R)$ back into itself it admits a unique fixed point in that set.
\end{Proof}

\begin{Proof}[Proof of Theorem \ref{thm-main}]
The existence of the undulated bilayer solution is a direct consequence of Proposition\,\ref{p:contraction}. It remains to establish its asymptotic form. 
Denoting the unscaled fixed point by $\phi_0^*:=\delta\veps^{1/4}\widetilde{\phi}_0^*$, we conclude from \eqref{def-mbi}, \eqref{e:undulation-sch}, \eqref{v-decomp}, and \eqref{e:scaling} that the undulated solution, $u_n$, takes the form
\begin{equation}\label{e:un-0}
    u_n=u_{h}\left(\frac{x_2}{\veps};\veps,\eta_d(x_1,\delta)\right)+\psi_0\phi_0^*+H(\phi_0^*;\delta,\veps),
\end{equation}
where from \eqref{e:Hexp} and \eqref{e:Hsch-ex} we have the estimate 
\begin{equation}\label{e:H*}
    \|H(\phi_0^*;\delta,\veps)\|_{H^4(\R^2)}=\caO(|\delta\veps|+|\veps|\|\phi_0^*\|_{H^4(\R)}+\|\phi_0^*\|_{H^4(\R)}^2).
\end{equation}
It remains to identify the leading order form of $\widetilde{\phi}_0^*$ and quantify the size of the remainder terms. We apply Lemma \ref{l:G*f}, together with the estimates \eqref{e:oK0-0}, \eqref{e:oK0-1},  and \eqref{e:tildeM}, to the right-hand side of the rescaled fixed point equality \eqref{e:tildeT} yielding, 
\begin{equation}\label{e:phi-est-fp}
    \|\widetilde{\phi}_0^*\|_{H^4(\R)}=\caO(1+|\veps|^{3/4}+|\veps|\|\widetilde{\phi}_0^*\|_{H^4}+|\delta|\,\|\widetilde{\phi}_0^*\|_{H^4}+|\delta\veps^{-3/4}|\,\|\widetilde{\phi}_0^*\|_{H^4}^2).
\end{equation}
From Young's inequality we deduce that the fixed point is $O(1)$. Returning to \eqref{e:tildeT} we absorb factors of $\|\widetilde{\phi}_0^*\|_{H^4(\R)}$ and eliminate $\delta=\widetilde{\delta}\veps^q$  to conclude that
\begin{equation}\label{e:phi-est}
    \widetilde{\phi}_0^*=-\veps^{3/4}\lambda_0^{-2}G*\oK_0+\caO(|\veps|^{3/4}+\widetilde{\delta}\veps^{q-3/4}).
\end{equation}
Combining \eqref{e:un-0}, \eqref{e:H*}, and \eqref{e:phi-est}, we have
\begin{equation}\label{e:un-1}
    u_n=u_{h}\left(\frac{x_2}{\veps};\veps,\eta_d(x_1,\delta)\right)-\delta\veps\lambda_0^{-2}\psi_0\Big(G*\oK_0\Big)+\caO(\delta\veps+\delta^2\veps^{-1/2}).
\end{equation}

Taking advantage of similar arguments as in Lemma \ref{l:G*f}, together with the expansions $A=1+\caO(\veps)$ and $B=\sqrt{-\alpha_0\veps}+\caO(\veps^{3/2})$ from \eqref{e:AB}, the leading order term $G*\oK_0$ can be evaluated directly. More specifically, we have
\[
\begin{aligned}
G*\oK_0(t)=&\int_\R G(t-s)\oK_0(s)\rmd s\\
=&\frac{1}
{4AB(A^2+B^2)}\int_\R\Big[ AE_e\big(B(t-s)\big)\cos\big(A(t-s)\big)+BE_o\big(B(t-s)\big)\sin\big(A(t-s)\big)\Big]\oK_0(s)\rmd s,\\
=&\frac{1}{4\sqrt{-\alpha_0\veps}}\int_\R E_e(B(t-s))\cos(A(t-s))\oK_0(s)\rmd s+\caO(\veps^{-1/4}).
\end{aligned}
\]
Using the double angle formula to break the $\cos$ term into $s$ and $t$ dependent parts gives the expression
\[
\begin{aligned}
G*\oK_0(t)
=&\frac{1} 
{4\sqrt{-\alpha_0\veps}}\bigg[\Big(\int_\R E_e\big(B(t-s)\big)\cos(As) \oK_0(s) \rmd s\Big)\cos(At) \\&\hspace{0.5in}+\Big(\int_\R E_e\big(B(t-s)\big)\sin(As) \oK_0(s) \rmd s\Big)\sin(At)\,\, \bigg]+\caO(\veps^{-1/4}),
\end{aligned}
\]
where the error estimate is in the $H^4(\R^2)$-norm.
The term $E_e$ is slowly varying since $B\ll1.$ Since $\oK_0$ has $\caO(1)$ compact support, localized near $0$,  we may approximate $E_e$ by its value at $s=0$, which affords the simplification
\beq
\label{e:Gf-exact}
G*\oK_0(t) = \frac{E_e(\sqrt{-\alpha_0\veps}t)}
{4\sqrt{-\alpha_0\veps}} \left[ \Big(\int_\R \oK_0(s)\cos(s)\rmd s\Big)\cos(At)+\Big(\int_\R \oK_0(s)\sin(s)\rmd s\Big)\sin(At)\right]+\caO(\veps^{-1/4}).
\eeq
From the form \eqref{K0-exp} of $\oK_0$ we may simplify, 
\[ 
\begin{aligned}
\int_\R  \oK_0(s)\cos(s)\rmd s = & \langle\psi_0,W'(u_0)\rangle_{L^2(\R)} \int_\R \left(\xi^{(4)}(s)+ 2\xi''(s)\right) \cos(s)\rmd s=-\langle\psi_0,W'(u_0)\rangle_{L^2(\R)}\Xi_{e,1},\\
\int_\R  \oK_0(s)\sin(s)\rmd s = & \langle\psi_0,W'(u_0)\rangle_{L^2(\R)} \int_\R \left(\xi^{(4)}(s)+ 2\xi''(s)\right) \sin(s)\rmd s=-\langle\psi_0,W'(u_0)\rangle_{L^2(\R)}\Xi_{e,1},
\end{aligned}
\]
where the Fourier coefficients $\Xi_{e/o,1}$ are defined in \eqref{e:xi-Fourier}.

Returning these results to \eqref{e:Gf-exact}, we conclude that
\beq
\label{e:Gf-exact2}
G*\oK_0(t)=\frac{\langle\psi_0,W^\prime(u_0) \rangle_{L^2(\R)}}{4\sqrt{-\alpha_0\veps}}\rme^{-\sqrt{-\alpha_0\veps}|t|}\bigg[\Xi_{\rme, 1}\cos(At)+\Xi_{\mathrm{o},1}\sin(At)
\bigg]+\caO(\veps^{-1/4}),
\eeq
where the error estimate is in the $H^4(\R^2)$-norm.  
Plugging \eqref{e:Gf-exact2} into \eqref{e:un-1}, we obtain \eqref{e:main-result}, which concludes the proof.
\end{Proof}

\subsection*{Acknowledgment}
K. Promislow acknowledges support from NSF-DMS grant 1813203. 
Q. Wu acknowledges support from  NSF-DMS grant 1815079.

\bibliographystyle{siam}
\bibliography{myref}

\end{document}